\documentclass[12pt,reqno]{amsart}
\usepackage{latexsym,amssymb,amsmath,multicol,rotating,
lscape, color}
 2

\usepackage{longtable}
\newtheorem{thm}{Theorem}[section]

\newtheorem{prop}[thm]{Proposition}
\newtheorem{cor}[thm]{Corollary}
\newcommand{\thmref}[1]{Theorem~\ref{#1}}

\newcommand{\rmkref}[1]{Remark~\ref{#1}}

\newcommand{\corref}[1]{Corollary~\ref{#1}}

\theoremstyle{remark}
\newtheorem{rmk}{Remark}[section]

\begin{document}

\title[Represenatations of squares]
{Representations of squares by certain diagonal quadratic forms in odd number of variables}

\author{B. Ramakrishnan, Brundaban Sahu and Anup Kumar Singh}
\address[B. Ramakrishnan]{Indian Statistical Institute, North-East Centre, Punioni, Solmara\\ 
Tezpur - 784 501, Assam, India} 
\address[Anup Kumar Singh and Brundaban Sahu]
{School of Mathematical Sciences, National Institute of Science 
Education and Research Bhubaneswar, HBNI, 
Via Jatni, Khurda, Odisha - 752 050,
India.}

\email[B. Ramakrishnan]{ramki@isine.ac.in, b.ramki61@gmail.com}
\email[Brundaban Sahu]{brundaban.sahu@niser.ac.in}
\email[Anup Kumar Singh]{anupsinghmath@gmail.com, anupsinghmath@niser.ac.in}

\date{\today}

\begin{abstract}
In this paper, we consider the following diagonal quadratic forms
\begin{equation*}
a_1x_1^2 + a_2x_2^2 + \cdots + a_{\ell}x_{\ell}^2,
\end{equation*}
where $\ell\ge 5$ is an odd integer and  $a_i\ge 1$ are integers. By using the extended Shimura correspondence, we obtain explicit formulas for the number of representations of $|D|n^2$ by the above type of quadratic forms, where $D$ is either a square-free integer or a fundamental discriminant such that $(-1)^{(\ell-1)/2}D > 0$. We demonstrate our method with many examples, in particular,  we obtain all the formulas (when $\ell =5$) obtained in the work of Cooper-Lam-Ye \cite{Ye1} and all the representation formulas for $n^2$ obtained in  \cite{Ye2} when $n$ is even. The works of Cooper et. al make use of certain theta function identities combined with a method of Hurwitz to derive these formulas. It is to be noted that our method works in general with arbitrary coefficients $a_i$. As a consequence to some of our formulas, we obtain certain identities among the representation numbers and also some congruences involving Fourier coefficients of certain newforms of weights $6, 8$ and the divisor functions. 
\end{abstract}

\maketitle

\section{Introduction}

\smallskip

Let $r_k(n)$ denote the number of representations of a positive integer $n$ as a sum of $k$ integer squares. Finding formulas for $r_k(n)$ is a classical problem in number theory. A general formula for $r_k(n)$ is not known so far when $k$ is odd.  
However, some explicit formulas of $r_k(n)$ or $r_k(n^2)$ for $k =1,3,5,7,9,11$ and $13$ are known. Here we give some references \cite{{hirsch}, {cooper1}, {cooper2}, {hurwitz}, {lomadze}, {sandham}, {sandham1}}. For a comprehensive list 
we refer to \cite{cooper2}. 
In \cite{gun}, the first author in collaboration with Sanoli Gun obtained a general 
formula for $r_{2m+1}(n^2)$ ($2\le m \in {\mathbb N}$) by using the extended Shimura correspondence obtained by A. G. van Asch \cite{asch}. 
The problem becomes difficult if one considers diagonal quadratic forms with integer coefficients (and one of the coefficients 
is  greater than 1).  In this direction, formulas for the number of representations of $n^2$ by the quadratic forms $x^2 + b y^2 + c z^2$ for $b,c \in \{1,2,3\}$ were obtained by S. Cooper and Y. Lam in \cite{cooper-lam}. Similar works to find formulas for the number of representations of $n^2$ by sums of 5 squares and 7 squares with coefficients (i.e., by the quadratic forms $a_1x_1^2 + \cdots + a_\ell x_\ell^2$, with $\ell = 5,7$, $a_i\in 
\{1,2,3,6\}$, when $\ell =5$ and $a_i\in \{1,2,4\}$, when $\ell =7$) were carried out by S. Cooper, Y. Lam and D. Ye in \cite{{Ye1}, {Ye2}} and they refer these quadratic forms as quinary ($\ell =5$) and septenary ($\ell =7$) forms. In these works, 
Cooper-Lam-Ye used a method of Hurwitz as the main tool to get explicit formulas for the number of representations.  

The purpose of this paper is to use the Shimura correspondence to deduce a more general formula. More precisely, we find formulas for $r_\ell(a_1,\ldots, a_\ell; |D|n^2) = r_\ell({\bf a}; |D|n^2)$, the number of representations of $|D|n^2$ by the quadratic form $a_1x_1^2 + \cdots + a_\ell x_\ell^2$, where ${\bf a} = (a_1, \ldots, a_\ell)$, $a_i$'s are positive integers, 
$\ell\ge 5$ is an odd integer and $D$ is either a square-free integer $t$ or a fundamental discriminant (which depends on $\ell$, $N_{\bf a}$). Here $N_{\bf a}$ is the least common multiple of all the coefficients $a_i$, $1\le i\le \ell$ and we assume that $(-1)^{(\ell-1)/2} D > 0$.  This is achieved by observing the fact that the generating function corresponding to the quadratic form is a modular form of weight $\ell/2$ on $\Gamma_0(4N_{\bf a})$ with certain quadratic character depending on ${\bf a}$. 
So, we can apply the extended Shimura map ${\mathcal S}_{D}$ obtained by Jagathesan and Manickam 
\cite{{mm1}, {mm2}} and express it as a linear combination of basis elements of  modular forms of weight $\ell -1$ and level $2N_{\bf a}$. Now the required formula follows by comparing the Fourier coefficients and taking M\"obius inversion. Thus, our approach gives a more general formula for the number of representations.

We now give the details of results proved in this paper. In \S 2, we give some preliminaries and state our main theorem and in \S 3 we give a proof. In \S 4, explicit examples are discussed for the cases $\ell = 5, 7, 9$.
In \S 4.1 we consider the quinary forms ($\ell =5$) and  provide 54 examples by taking $a_i \in \{1,2,3,4,6\}$. 
All the 18 formulas obtained in the work of Cooper-Lam-Ye \cite{Ye1} are presented as \corref{quinary}. Necessary details for the remaining 36 formulas are given in Tables 4 and 5. Some explicit examples (among these 36 cases) are presented in 
Corollaries 4.2 and 4.3. The septenary case ($\ell =7$) is presented in \S 4.2. Here we take $a_i\in \{1,2,4\}$, which gives 27 forms, out of which 18 cases correspond to the formulas obtained in \cite{Ye2}. Our formulas coincide with the formulas of Cooper-Lam-Ye when $n$ is even. We remark that the existing results on the mapping property of Shimura maps allow us to use only the Shimura-Kohnen map when $N_{\bf a}$ is even.  However, by assuming the mapping property of the extended Shimura maps ${\mathcal S}_t$ in the case when $\ell \equiv 3\pmod{4}$ and $N_{\bf a}$ is even, one could obtain formulas for the number of representations of $n^2$. In particular, the remaining formulas obtained in \cite{Ye2} could be proved under this assumption. This is explained in \S 4.2.2 and the conjectural formulas are given in this section.  In \S 4.3, we consider the case $\ell =9$ and obtain 52 formulas by taking $a_i$ in the three sets $\{1,b\}$, $b=2,3$ and $\{1,2,4\}$. Complete data for obtaining the formulas are given in Tables 7 to 10 and sample formulas are provided in \corref{9var}. As a consequence to these formulas, we derive certain congruences between newform Fourier coefficients and the divisor function, which is presented in \corref{cong}. We also obtain certain relations among the representation numbers $r_9({\bf a};n^2)$, which is presented as a proposition in \S 4.3. So far the known results for the number of representations are mainly for perfect squares. 
The method of Shimura correspondence gives such formulas for the number of representations of integers which are not necessarily perfect squares and we could get formulas for $tn^2$, where $t$ is a positive square-free integers. 
In \S 5 we discuss about such formulas (when $\ell =7$) and list some of them which are not proved earlier. 
We also mention about additional formulas that can be obtained using our method which are not presented here due to the 
large size of the data. 

\smallskip

\section{Preliminaries and Main Results}

We consider the following quadratic form in $\ell$ variables with (positive) integer coefficients $a_i$:
\begin{equation}\label{1}
a_1x_1^2 + a_2x_2^2 + \cdots + a_{\ell}x_{\ell}^2.
\end{equation}
Let $r_{\ell}(a_1, a_2,\ldots, a_\ell; n)$ denote the number of representations of a natural number $n$ by the above quadratic form. i.e., writing ${\bf a} = (a_1,a_2,\ldots, a_\ell)$, 
\begin{equation}\label{eq1}
r_\ell({\bf a}; n) := r_{\ell}(a_1,a_2,\ldots, a_\ell; n)= \# \left\{(x_1,x_2,\ldots,x_\ell) \in \mathbb{Z}^\ell; \displaystyle\sum_{i=1}^{\ell}a_i x_i^2 = n \right\}.
\end{equation}
When all the $a_i$'s are equal to $1$, it is denoted by $r_\ell(n)$. If $a_j$ appears with multiplicity $i_j$, then we represent 
${\bf a}$ by  $(a_1^{i_1}, a_2^{i_2}, \cdots, a_{r}^{i_r})$, where  $i_1+i_2+\cdots +i_r= \ell$. We drop the power 
$i_j$, if it is equal to 1. 

\smallskip

Before stating our main theorem, we shall fix some notations. For natural numbers $k$ and $M$, the vector spaces of modular forms of weight $k+1/2$ (resp. weight $2k$) on $\Gamma_0(4M)$ with character $\chi$ modulo $4M$ (resp. 
on $\Gamma_0(M)$ with character $\chi$ modulo $M$) as $M_{k+1/2}(4M,\chi)$ (resp. $M_{2k}(M,\chi)$). 
The respective subspaces of cusp forms are denoted by $S_{k+1/2}(4M,\chi)$ and $S_{2k}(M,\chi)$. In the case of principal  character, we omit $\chi$ from the notation. 
For a vector ${\bf a} = (a_1, a_2, \ldots, a_\ell)$, let $N_{\bf a}$ be as defined in the introduction (lcm of all the integers $a_1, \ldots, a_\ell$). Corresponding to ${\bf a}$, let us also associate a quadratic Dirichlet character $\prod_{j=1}^\ell \left(\frac{a_j}{\cdot}\right)$ and denote it by $\psi_{\bf a}$. Note that $\psi_{\bf a}$ is the principal Dirichlet character modulo $4N_{\bf a}$, if for each $j$, $1\le j\le \ell$, either $a_j$ is a perfect square or it appears even number of times in \eqref{1}. 

\smallskip

\begin{thm}\label{thm:main}
For a vector ${\bf a} = (a_1, a_2, \ldots, a_\ell)$, with $a_i\in {\mathbb N}$, let $N_{\bf a}$ be the positive integer and 
$\psi_{\bf a}$ be the Dirichlet character associated to ${\bf a}$ as defined above. Let $\nu(\ell,N_{\bf a})$ denote the dimension of the vector space $M_{\ell-1}(2N_{\bf a})$. Let $D$ be a square-free integer or a fundamental discriminant depending on $N_{\bf a}$ is odd or even such that $(-1)^{(\ell-1)/2} D > 0$. Then for $n\in {\mathbb N}$, 
we have 
\begin{equation} \label{eq:main}
r_{\ell}({\bf a}; |D|n^2) =  
\sum_{d\vert n\atop{(d,2N_{\bf a})=1}}\mu(d) \psi_{\bf a}(d) \left(\frac{D}{d}\right) d^{(\ell-3)/2} \sum_{j=1}^{\nu(\ell,N_{\bf a})} 
\lambda_{\ell, D,j}({\bf a}) A_{\ell, {\bf a};j}(n/d),
\end{equation}
where $\lambda_{\ell, D,j}({\bf a})$ are constants which depend on $\ell, D, {\bf a}$ and also on the choice of a basis of 
$M_{\ell -1}(2N_{\bf a})$. The $n$-th Fourier coefficients of these basis elements are denoted by   $A_{\ell,{\bf a};j}(n)$, $1\le j\le \nu(\ell, N_{\bf a})$. 
\end{thm}

\smallskip

\section{Proof of the Main Theorem}
\smallskip

In \cite{{mm1},{mm2}}, Jagathesan and Manickam have obtained the Shimura correspondence for non-cusp forms of half-integral weight on $\Gamma_0(4N)$, $N\in {\mathbb N}$. When $N=1$, the work was carried out by A. G. van Asch \cite{asch}. For more details on Shimura and Shimura-Kohnen maps, we refer to the works of Shimura and Kohnen 
\cite{{shimura},{kohnen},{kohnen1}}. Below we mention a result of Jagathesan and Manickam, which is the main ingredient in our proof.

\begin{thm} {\rm (}Jagathesan and Manickam  \cite{{mm1}, {mm2}}{\rm )} \label{thm-mm}\\
Let $k\ge 2$ be an integer and $t$ be a square-free integer with $(-1)^k t>0$. Write $D=t$ or $4t$ 
according as $t\equiv 1$ or $2,3 \pmod{4}$. Then $D$ is a fundamental discriminant 
with $(-1)^k D>0$. For a modular form $f(z) = \sum_{n\ge 0} a(n) q^n \in M_{k+1/2}(4N,\chi)$, where $\chi$ is an even 
quadratic Dirichlet character modulo $4N$, define the $t$-th Shimura map and $D$-th Shimura-Kohnen map as follows:
\begin{equation}\label{shimura-map}
{\mathcal S}_t(f)(z) = a(0) H_k(t) + \sum_{n\ge 1} \left( \sum_{d\vert n\atop {(d,2N)=1}} 
\chi(d)  \left(\frac{t}{d}\right)d^{k -1} a(|t|n^2/d^2)\right) ~q^n,
\end{equation}
\begin{equation}\label{shimura-kohnen-map}
{\mathcal S}_D(f)(z) = a(0) H_k(D) + \sum_{n\ge 1} \left( \sum_{d\vert n\atop {(d,N)=1}}
\chi(d)   \left(\frac{D}{d}\right)d^{k -1} a(|D|n^2/d^2)\right) ~q^n,
\end{equation}
where $z\in {\mathbb C}$, ${\rm Im} z>0$, $q = e^{2\pi iz}$ and 

\begin{equation}
\begin{split}
H_k(t) &= \frac{1}{2} L\left(1-k, \left(\frac{D}{\cdot}\right)\right) 
\prod_{p\vert 2N}\left(1 - \left(\frac{D}{p}\right)p^{k -1}\right), \\
H_k(D) &= \frac{1}{2} L\left(1-k, \left(\frac{D}{\cdot}\right)\right) 
\prod_{p\vert N}\left(1 - \left(\frac{D}{p}\right)p^{k -1}\right).
\end{split}
\end{equation} 
Then, we have ${\mathcal S}_t(f) \in M_{2k}(2N)$, if $N$ is odd and 
${\mathcal S}_D(f) \in M_{2k}(2N)$, if $N$ is even.
\end{thm}

\begin{rmk}
In the above theorem, extended Shimura map is obtained for the modular forms space with trivial character. However, 
Professor Manickam informed us that in his forthcoming work the results are proved for the spaces with real quadratic character. So, we make use of this fact also in our work. 
\end{rmk}

\bigskip

\noindent {\bf Proof of \thmref{thm:main}}. 

We are now ready to prove our main theorem. 
Let $\ell \ge 5$ be an odd integer.  Let $\Theta(z) = \sum_{n\in {\mathbb Z}} q^{n^2}$ be the classical theta series of weight $1/2$ for the group $\Gamma_0(4)$. Then the generating function for the quadratic form 
$a_1x_1^2 + \cdots + a_\ell x_\ell^2$ is the following product of theta functions:
$$
\Theta_{\bf {a}}(z) = \prod_{j=1}^\ell \Theta(a_j z),
$$
where ${\bf a}  = (a_1,a_2,\ldots, a_\ell)$ and it is a  modular form of weight $\ell/2$ on $\Gamma_0(4N_{\bf a})$, 
with character $\psi_{\bf a}$ where $N_{\bf a}$ and $\psi_{\bf a}$ are defined as before. This follows from the basic theory of modular forms. For a proof of this fact, we refer to \cite[Fact II, p.758]{rss1}. Now we make use of \thmref{thm-mm}. 
When $N_{\bf a}$ is odd, we apply the Shimura map 
${\mathcal S}_t$ with square-free $t$ such that $(-1)^{(\ell-1)/2}t > 0$ on $\Theta_{\bf a}(z)$ and when $N_{\bf a}$ is even, 
then we apply the Shimura-Kohnen map ${\mathcal S}_{D}$, where $D$ is a fundamental discriminant such the $(-1)^{(\ell-1)/2}D>0$ on the theta function $\Theta_{\bf a}(z)$. For the sake of uniformity in the theorem, we used the symbol $D$ to denote a square-free integer or a fundamental discriminant depending on the parity (odd or even) of $N_{\bf a}$. Thus, ${\mathcal S}_D(\Theta_{\bf a})(z)$ is a modular form in $M_{\ell-1}(2N_{\bf a})$. Let $f_{\ell,j}(z)$, $1\le j\le \nu(\ell, N_{\bf a})$ be a basis for the space $M_{\ell-1}(2N_{\bf a})$, whose $n$-th Fourier coefficients are denoted by $A_{\ell, {\bf a}; j}(n)$, $1\le j\le 
\nu(\ell, N_{\bf a})$. Therefore, we can write ${\mathcal S}_D(\Theta_{\bf a})(z)$ as a linear combination of these basis elements and then comparing their $n$-th Fourier coefficients, we get the following identity for any integer $n\ge 1$:
\begin{equation}
\sum_{d\vert n\atop {(d,2N_{\bf a})=1}} \psi_{\bf a}(d)  \left(\frac{D}{d}\right) d^{(\ell -3)/2} r_{\ell}({\bf a}; |D|n^2/d^2) 
= \sum_{j=1}^{\nu(\ell, N_{\bf a})} \lambda_{\ell, D, j}({\bf a}) A_{\ell, {\bf a};j}(n),
\end{equation}
where $\lambda_{\ell, D, j}({\bf a})$ are constants that depend on $\ell$, $D$, ${\bf a}$ and the basis chosen for the space $M_{\ell-1}(2N_{\bf a})$. Taking M\"obius inversion, the above identity gives the required expression for $r_\ell({\bf a};|D|n^2)$ as stated in the theorem, namely:
\begin{equation}
r_\ell({\bf a}; |D|n^2) = \sum_{d\vert n\atop {(d,2N_{\bf a})=1}} \mu(d) \psi_{\bf a}(d)  \left(\frac{D}{d}\right) d^{(\ell -3)/2} 
\sum_{j=1}^{\nu(\ell, N_{\bf a})} \lambda_{\ell, D, j}({\bf a}) A_{\ell, {\bf a};j}(n/d).
\end{equation}
This completes the proof. 

\bigskip

In the following sections, we make the above formulas more explicit. 

\section{Sample formulas} 

In this section we shall apply our main theorem to derive explicit formulas for $r_{\ell}({\bf a};n^2)$ for $\ell =5, 7, 9$, with certain coefficients $a_i$. In the case of quinary forms ($\ell =5$), we prove all the identities obtained in \cite{Ye1} and in the case of septenary forms ($\ell =7$), we prove most of the formulas obtained in \cite{Ye2} apart from proving new identities. 
We add a section which contains some conjectural formulas, which are derived under the assumption of the mapping property of the $t$-th Shimura maps.   When $\ell =9$, we get the formulas when $a_i $ belong to the sets $\{1,2\}$, $\{1,3\}, \{1,2,4\}$. 
All the formulas obtained in our work are listed as Corollaries to our main theorem. 

\smallskip

We use the following notations in our formulas. For an odd positive integer $\displaystyle{m=\prod_{p\ge 3} p^{\lambda_p}}$, we define the following functions. 

\begin{align}
s_1(m) & = \prod_{p\ge 3} \left(\frac{p^{3\lambda_p+3}-1}{p^3-1} - p \frac{p^{3\lambda_p}-1}{p^3-1}\right),\label{s1}\\
s_2(m) & = \prod_{p\ge 3} \left(\frac{p^{3\lambda_p+3}-1}{p^3-1} - p \left(\frac{2}{p}\right) \frac{p^{3\lambda_p}-1}{p^3-1}\right),
\label{s2}\\
s_3(m) & = \prod_{p\ge 3} \left(\frac{p^{3\lambda_p+3}-1}{p^3-1} - p \left(\frac{3}{p}\right) \frac{p^{3\lambda_p}-1}{p^3-1}\right),
\label{s3}\\
s_4(m) & = \prod_{p\ge 3} \left(\frac{p^{3\lambda_p+3}-1}{p^3-1} - p \left(\frac{6}{p}\right) \frac{p^{3\lambda_p}-1}{p^3-1}\right),
\label{s4}
\end{align}

\begin{align}
s_5(m) &= \prod_{p\ge 3} \left(\frac{p^{5\lambda_p+5} -1}{p^5-1} - p^2 \left(\frac{-2}{p}\right) 
\frac{p^{5\lambda_p}-1}{p^5-1}\right),\label{sm}\\
s_6(m) &= \prod_{p\ge 3} \left(\frac{p^{5\lambda_p+5} -1}{p^5-1} - p^2 \left(\frac{-4}{p}\right) 
\frac{p^{5\lambda_p}-1}{p^5-1}\right).\label{tm}\\
s_7(m) & = \prod_{p\ge 3} \left(\frac{p^{7\lambda_p+7}-1}{p^7-1} - p^3 \frac{p^{7\lambda_p}-1}{p^7-1}\right).\label{s7}
\end{align}

\begin{rmk}
The functions \eqref{s1}-\eqref{s4} appear in the formulas of Theorems 1.1 to 1.4 of \cite{Ye1}. 
We note that the functions \eqref{sm} and \eqref{tm} are denoted as $s(m)$ and $t(m)$ respectively in  \cite[Eqs. (2.1), (2.2)]{Ye2}.
\end{rmk}

\subsection{Quinary forms}

In this section, we assume that $\ell =5$ and consider the form $a_1 x_1^2 + \cdots + a_5 x_5^2$, where $a_i$ belong to the set $\{1,2,3,4,6\}$. Specifically, we take them in the two sets $\{1,2,4\}$ and $\{1,2,3,6\}$. There are 14 forms corresponding to the first set and the functions corresponding to them belong to the space of modular forms of weight $5/2$ 
on $\Gamma_0(8)$ (if $a_i \not= 4$) or $\Gamma_0(16)$ (if some $a_i =4$) with character depending on ${\bf a}$. There are 40 quinary forms when $a_i$'s belong to the set $\{1,2,3,6\}$. The corresponding generating functions belong to either $M_{5/2}(12,\chi)$, $M_{5/2}(16,\chi')$ or $M_{5/2}(24,\chi'')$, where $\chi$, $\chi'$, $\chi''$ are real quadratic characters modulo $12$, $16$ or $24$ depending on the coefficients $a_i$'s. To apply our method, we need to find explicit bases for the following modular form spaces: $M_4(N)$, where $N=4,6,8,12$. In this section, we consider certain bases for these modular form spaces and apply our theorem to get formulas for $r_5({\bf a};n^2)$ for the above $54$ quinary forms. In particular, we derive all the 18 formulas proved in \cite{Ye1} and the remaining 36 formulas are new. 
In the table below, we shall give the list of vectors and the corresponding modular forms space. In Table 2, we shall describe the required basis elements for the evaluation of our formulas. As mentioned earlier, we use the notation ${\bf a} = (a_1^{i_1}, a_2^{i_2}, \cdots, a_{r}^{i_r})$, $i_1+i_2+\cdots +i_r= \ell$.

\smallskip

\begin{center}
\centerline{{\bf {Table 1}} ($\ell=5$, $a_i \in \{1,2,4\}, \{1,2,3,6\}$)}
{\Small
\begin{longtable}{|c|l|}
\hline 
${\bf a} = (a_1^{i_1}, a_2^{i_2}, \cdots, a_{r}^{i_r})$ with $i_1+i_2+\cdots +i_r= 5 $  &$4N_{\bf a}, \psi_{\bf a}$  \\
\hline
$(1^3,2^2), (1,2^4)$ & 8, $\chi_0$ \\
$(1^4,2), (1^2,2^3)$ & 8, $\chi_2$ \\
\hline
$ (1^4,4), (1^3,4^2)$, $(1^2,4^3), (1,4^4), (1,2^2,4^2), (1^2,2^2,4) $ & 16, $\chi_0$\\
$ (1^3,2,4), (1^2,2,4^2), (1,2^3,4) , (1,2,4^3)$, & 16, $\chi_2$\\
\hline
$(1,3^4), (1^3,3^2)$ & 12, $\chi_0$\\
$(1^4,3), (1^2,3^3)$ & 12, $\chi_3$\\
\hline 
$(1,2^2,3^2), (1,3^2,6^2), (1^3,6^2), (1,6^4),(1,2^2,6^2), (2^3,3,6), (2,3^3,6), (2,3,6^3), (1^2,2,3,6)$
& 24, $\chi_0$\\  
 $(2^3,3^2), (2,3^4), (1^2,2,3^2), (1^3,3,6), (1,3^3,6), (1,3,6^3), (1^2,2,6^2), (2,3^2,6^2), (1,2^2,3,6)$
 &24, $\chi_2$\\  
 $(2^4,3), (2^2,3^3), (1^2,2^2,3),(1^2,3,6^2), (1^3,2,6), (1,2^3,6), (1,2,6^3), (2^2,3,6^2), (1,2,3^2,6)$ 
 &24, $\chi_3$\\  
 $(1^3,2,3), (1,2^3,3), (1,2,3^3), (1^4,6), (1^2,6^3),  (1^2, 3^2,6),(1^2,2^2,6),  (2^2,3^2,6),   (1,2,3,6^2)$
 & 24, $\chi_6$\\
\hline
\end{longtable}
}
\end{center}

\noindent In the above, $\chi_0$ denotes the principal character modulo $4N_{\bf a}$ and $\chi_m$ is the quadratic Dirichlet character $\left(\frac{m}{\cdot}\right)$ modulo $4m$, where $m\vert N_{\bf a}$. 

\bigskip

In the following table, we describe the basis elements for the vector spaces $M_{4}(N)$, $N = 4, 6, 8$ and $12$. 
We use the following notation for these basis elements. The function $E_k(z)$ denotes the normalised Eisenstein series 
of an even integer weight $k\ge 4$ for the full modular group whose Fourier expansion is given by 
\begin{equation}
E_k(z) = 1 - \frac{2k}{B_k} \sum_{n\ge 1} \sigma_{k-1}(n) q^n,
\end{equation} 
where $B_k$ is the $k$-th Bernoulli number, $q = e^{2\pi iz}$ and $\sigma_{k-1}(n)$ is the divisor function. 
For a given weight $k\ge 2$ and level $N$, if the normalised newform is unique, then we denote it by  $\Delta_{k,N}(z)$, and denote by $\tau_{k,N}(n)$ its $n$-th Fourier coefficient. If there are more than one newform of weight $k$, level $N$, then we denote them by $\Delta_{k,N;j}(z)$, where $j$ runs over the number of such newforms. Their $n$-th Fourier coefficients are denoted by $\tau_{k,N;j}(n)$. \\

\begin{center}
\centerline{{\bf Table 2}}

\bigskip

\begin{tabular}{|c|c|c|}
\hline 
 Level&  & Dimension\\
$2N_{\bf a}$&Basis elements &$\nu(5,N_{\bf a})$\\
\hline 
$4$ & $E_4(z),E_4(2z),E_4(4z)$ & 3\\
\hline 
$6$ & $E_4(z), E_4(2z), E_4(3z), E_4(6z), \Delta_{4,6}(z)$ & 5\\
\hline 
$8$ & $ E_4(z), E_4(2z), E_4(4z), E_4(8z), \Delta_{4,8}(z) $ & 5 \\
\hline 
$12$ &$ E_4(tz);t|12, \Delta_{4,6}(z), \Delta_{4,6}(2z), \Delta_{4,12}(z) $ & 9\\
\hline
\end{tabular}
\end{center}

\smallskip

\bigskip

Let $\eta(z) = q^{1/24} \displaystyle{\prod_{n=1}^\infty} (1-q^n)$ be the Dedekind eta-function, which is a modular form of weight $1/2$. Then the newforms $\Delta_{k,N}(z)$ appearing in the above table are given by the following eta-products or 
linear combination of eta-quotients. 
\begin{equation*}
\begin{split}
\Delta_{4,6}(z) & = \eta^2(z)\eta^2(2z)\eta^2(3z)\eta^2(6z); \qquad \quad 
\Delta_{4,8}(z) ~ =~ \eta^4(2z)\eta^4(4z);\\
\Delta_{4,12}(z)& =  \frac{\eta^2(2z)\eta^3(3z)\eta^3(4z)\eta^2(6z)}{\eta(z)\eta(12z)} - \frac{\eta^3(z) \eta^2(2z)\eta^2(6z)\eta^3(12z)}{\eta(3z)\eta(4z)}.
\end{split}
\end{equation*}

\smallskip

\bigskip

In the following, we derive the 18 formulas for the quinary forms obtained in the work of Cooper-Lam-Ye \cite{Ye1}. 
As $\ell =5$, $(\ell -1)/2$ is even, we may take $t=D =1$. In this case ${\mathcal S}_t$ and ${\mathcal S}_D$ are the same 
map and the image belongs to one of the spaces $M_{4}(N)$, where $N=4,6,8$ or $12$. By considering the 18 vectors ${\bf a}$ and taking $t=D=1$ in \eqref{eq:main} we express $r_5({\bf a};n^2)$ explicitly. In the following table, we first give the images 
of $\Theta_{\bf a}(z)$ under ${\mathcal S}_1$ with respect to the bases given in the previous table for ${\bf a} \in \{(1^4,2), (1^4,3), (1^3,2,3), (1^3,2^2), (1^4,4), (1^3,4^2), (1^2,2^2,4), (1^2, 4^3), (1,2^4), (1,2^2,4^2), 
(1,4^4), 
\linebreak 
(1,2^2,3^2),  (1,3^4), (1^2,2^3), (1,2^3,4), (1^2,2,3^2), (1,3^3,6), (1,2^3,6)\}$. 
Using the explicit images given in the table below (Table 3) and comparing the $n$-th Fourier coefficients and finally taking M\"obius inversion we obtain the 18 formulas corresponding to the vectors ${\bf a}$ listed above. 
We present these formulas in the form of a corollary after the table. \\

\bigskip
\bigskip
\vfill 

\centerline{{\bf Table 3}}

\begin{center}
\vspace{0.15cm}
{\small 
\begin{tabular}{|c|c|}
\hline 
${\bf a}$ & ${\mathcal S}_1(\Theta_{\bf a})(z)$\\
\hline 
$(1^4,2)$ & $\frac{1}{30} E_4(z) - \frac{8}{15} E_4(4z)$\\
\hline
$(1^4,3)$ & $\frac{1}{30} E_4(z) - \frac{2}{15} E_4(2z)  + \frac{3}{10} E_4(3z) - \frac{6}{5} E_4(6z) - \frac{18}{5} E_4(12z)$\\
\hline
$(1^3,2,3)$ & $\frac{1}{40} E_4(z) - \frac{1}{20} E_4(2z)  - \frac{9}{40} E_4(3z) + \frac{2}{5} E_4(4z) + \frac{9}{20} E_4(6z)$\\
\hline
$(1^3,2^2)$ & $\frac{1}{40} E_4(z) - \frac{1}{60} E_4(2z)$\\
\hline
$(1^4,4)$ & $\frac{1}{30} E_4(z) - \frac{23}{120} E_4(2z) + \frac{1}{5} E_4(4z)$\\
\hline
$(1^3,4^2)$ & $\frac{1}{40} E_4(z) - \frac{11}{60} E_4(2z) + \frac{1}{5} E_4(4z)$\\
\hline
$(1^2,2^2,4)$ & $\frac{1}{60} E_4(z) - \frac{1}{24} E_4(2z) + \frac{1}{15} E_4(4z)$\\
\hline
$(1^2,4^3)$ & $\frac{1}{60} E_4(z) - \frac{13}{120} E_4(2z) + \frac{2}{15} E_4(4z)$\\
\hline
$(1,2^4)$ & $\frac{1}{120} E_4(z) + \frac{1}{30} E_4(2z)$\\
\hline
$(1,2^2,4^2)$ & $\frac{1}{120} E_4(z) - \frac{1}{30} E_4(2z) + \frac{1}{15} E_4(4z)$\\
\hline
$(1,4^4)$ & $\frac{1}{120} E_4(z) - \frac{1}{60} E_4(2z) + \frac{1}{15} E_4(4z)$\\
\hline
$(1,2^2,3^2)$ & $\frac{1}{120} E_4(z) - \frac{1}{30} E_4(2z)  + \frac{3}{40} E_4(3z) - \frac{3}{20} E_4(6z)$\\
\hline
$  (1,3^4)$ & $\frac{1}{120} E_4(z) - \frac{11}{120} E_4(3z)$\\
\hline
$(1^2,2^3)$& $\frac{1}{60} E_4(z) +\frac{1}{60} E_4(2z) - \frac{16}{30} E_4(4z)$\\
\hline
$(1,2^3,4)$ & $\frac{1}{120} E_4(z) - \frac{1}{120} E_4(2z)  + \frac{1}{30} E_4(4z) - \frac{16}{30} E_4(8z)$\\
\hline
$(1^2,2,3^2)$ & $\frac{1}{60} E_4(z) - \frac{1}{30} E_4(2z)  - \frac{3}{20} E_4(3z) + \frac{4}{15} E_4(4z) + \frac{3}{10} E_4(6z)
-\frac{12}{5} E_4(12z)$\\
\hline
$(1,3^3,6)$ & $\frac{1}{120} E_4(z) - \frac{1}{60} E_4(2z)  - \frac{17}{120} E_4(3z) + \frac{2}{15} E_4(4z) + \frac{17}{60} E_4(6z) -\frac{34}{15} E_4(12z)$\\
\hline
$(1,2^3,6)$ & $\frac{1}{120} E_4(z) - \frac{1}{60} E_4(2z)  - \frac{3}{40} E_4(3z) + \frac{2}{15} E_4(4z) + \frac{3}{20} E_4(6z) 
-\frac{6}{5} E_4(12z)$\\
\hline
\end{tabular}
}
\end{center}

\smallskip

\begin{cor}\label{quinary}
For a natural number $n$, we have 
\begin{equation*}
\begin{split}
r_5(1^4, 2;n^2) & = \sum_{d \vert n} \mu(d)d\left( \frac{2}{d}\right) \left\{ 8 \sigma_3(n/d) - 128 \sigma_3(n/4d) \right\},\\
r_5(1^4, 3;n^2) & = \sum_{d \vert n} \mu(d)d\left( \frac{3}{d}\right) \left\{ 8 \sigma_3(n/d)-32 \sigma_3(n/2d) +72 \sigma_3(n/3d)- 288 \sigma_3(n/6d) \right\},\\
r_5(1^3, 2, 3;n^2) & = \sum_{d \vert n} \mu(d)d\left( \frac{6}{d}\right) \left\{ 6 \sigma_3(n/d)-12 \sigma_3(n/2d) -54 \sigma_3(n/3d)+96 \sigma_3(n/4d)\right.\\
& \left. \qquad +108  \sigma_3(n/6d)-864 \sigma_3(n/12d) \right\},\\
r_{5}(1^3,2^2 ; n^2)& =  \sum _{d|n}\mu(d)d\left(\frac{4}{d}\right)\left\{6\sigma_3(n/d)-4\sigma_3(n/2d)\right\},\\
r_{5}(1^4,4 ; n^2)&= \sum _{d|n}\mu(d)d\left(\frac{4}{d}\right)\left\{8\sigma_3(n/d)-46\sigma_3(n/2d)+ 48\sigma_3(n/4d)\right\},\\
r_{5}(1^3,4^2 ; n^2)& =  \sum_{d|n}\mu(d)d\left(\frac{4}{d}\right)\left\{6\sigma_3(n/d)-44\sigma_3(n/2d)+ 48\sigma_3(n/4d)
\right\},\\
r_{5}(1^2, 2^2, 4 ; n^2)&= \sum _{d|n}\mu(d)d\left(\frac{4}{d}\right)\left\{4\sigma_3(n/d)-10\sigma_3(n/2d)+ 16\sigma_3(n/4d)\right\},
\end{split}
\end{equation*}
\begin{equation*}
\begin{split}
r_{5}(1^2,4^3 ; n^2)&= \sum _{d|n}\mu(d)d\left(\frac{4}{d}\right)\left\{4\sigma_3(n/d)-26\sigma_3(n/2d)+ 32\sigma_3(n/4d)\right\},
\\
r_{5}(1,2^4 ; n^2)&= \sum _{d|n}\mu(d)d\left(\frac{4}{d}\right)\left\{2\sigma_3(n/d)+8\sigma_3(n/2d)\right\},\\
r_{5}(1, 2^2, 4^2 ; n^2)&= \sum _{d|n}\mu(d)d\left(\frac{4}{d}\right)\left\{2\sigma_3(n/d)-8\sigma_3(n/2d)+ 16\sigma_3(n/4d)\right\},\\
r_{5}(1,4^4 ; n^2)&= \sum _{d|n}\mu(d)d\left(\frac{4}{d}\right)\left\{2\sigma_3(n/d)-4\sigma_3(n/2d)+ 16\sigma_3(n/4d)\right\},\\
r_{5}(1,2^2,3^2 ; n^2)&= \sum_{d|n}\mu(d)d\left(\frac{4}{d}\right)\left\{2\sigma_3(n/d) - 4 \sigma_3(n/2d)+18 \sigma_3(n/3d) -36 \sigma_3(n/6d)\right\}, \\
r_5(1, 3^4;n^2) & = \sum_{d \vert n} \mu(d)d\left( \frac{4}{d}\right) \left\{ 2 \sigma_3(n/d)-22 \sigma_3(n/3d) \right\},\\
r_5(1^2, 2^3;n^2) & = \sum_{d \vert n} \mu(d)d\left( \frac{2}{d}\right) \left\{ 4 \sigma_3(n/d)+ 4 \sigma_3(n/2d) - 128 \sigma_3(n/4d) \right\},\\
r_5(1, 2^3, 4;n^2) & = \sum_{d \vert n} \mu(d)d\left( \frac{2}{d}\right) \left\{ 2 \sigma_3(n/d)-2 \sigma_3(n/2d) +8 \sigma_3(n/4d)- 128 \sigma_3(n/8d) \right\},\\
r_5(1^2, 2, 3^2;n^2) & = \sum_{d \vert n} \mu(d)d\left( \frac{2}{d}\right) \left\{ 4 \sigma_3(n/d)-8 \sigma_3(n/2d) -36 \sigma_3(n/3d)+64 \sigma_3(n/4d)\right.\\
& \left. \qquad +72 \sigma_3(n/6d)-576 \sigma_3(n/12d) \right\},\\
r_5(1, 3^3, 6;n^2) & = \sum_{d \vert n} \mu(d)d\left( \frac{2}{d}\right)\{2 \sigma_3(n/d) - 4 \sigma_3(n/2d) - 34 \sigma_3(n/3d)
+ 32 \sigma_3(n/4d) \\
& \qquad + 68 \sigma_3(n/6d) - 544 \sigma_3(n/12d)\},\\
r_5(1, 2^3, 6;n^2) & = \sum_{d \vert n} \mu(d)d\left( \frac{3}{d}\right)\{ 2 \sigma_3(n/d) -4 \sigma_3(n/2d) - 18 \sigma_3(n/3d) + 32 \sigma_3(n/4d)\\
& \qquad  +36 \sigma_3(n/6d) - 288 \sigma_3(n/12d)\}.\\
\end{split}
\end{equation*}
\end{cor}

\smallskip

\begin{rmk}
All the above 18 formulas give rise to the same formulas obtained in the work of Cooper-Lam-Ye \cite{Ye1}. We demonstrate this in the next section.
\end{rmk}

\smallskip

\subsubsection{\bf Equivalence of formulas}

\smallskip

\noindent 
As mentioned in the above remark, all the above formulas are the same as given in \cite{Ye1}. We will demonstrate the equivalence in three cases. The rest of them follow using similar computations. When ${\bf a} = (1^4,2)$, we have 


\begin{equation*}
r_5(1^4, 2;n^2)  = \sum_{d \vert n} \mu(d)d\left( \frac{2}{d}\right) \left\{ 8 \sigma_3(n/d) - 128 \sigma_3(n/4d) \right\}.
\end{equation*}
Writing $n= 2^{\lambda_2} m$, $m =\displaystyle{\prod_{p\vert m\atop{p\ge 3}}} p^{\lambda_p}$, 
the above formula becomes, when $\lambda_2\ge 1$,
\begin{equation*}
\begin{split}
r_5(1^4, 2;n^2)  &= \sum_{d \vert n} \mu(d)d\left( \frac{2}{d}\right) \left\{ 8 \sigma_3(n/d) - 128 \sigma_3(n/4d) \right\}\\
& = \big(8 \sigma_3(2^\lambda) - 128 \sigma_2(2^{\lambda -2})\big) \sum_{d\vert m} \mu(d) \left(\frac{2}{d}\right) d \sigma_3(m/d),\\
& = 24 \left(\frac{2^{3\lambda_2 +1} + 5}{2^3-1}\right) s_2(m),
\end{split}
\end{equation*}
where $s_2(m)$ is given by \eqref{s2}. When $\lambda_2 =0$, i.e., when $n$ is odd, we only get the factor $8$ (outside the product) in the last equation. This is exactly the formula obtained in \cite[Theorem 1.2]{Ye1}.  Next we consider the case ${\bf a} = (1^4,,3)$. Writing $n= 2^{\lambda_2} 3^{\lambda_3} m$, with $\gcd(m,6)=1$,  we have 

\begin{equation*}
\begin{split}
r_5(1^4, 3;n^2) & = \sum_{d \vert n} \mu(d) d \left(\frac{3}{d}\right) \left\{8 \sigma_3(n/d)-32 \sigma_3(n/2d) +72 
\sigma_3(n/3d)- 288 \sigma_3(n/6d) \right\}\\
& = 8 \big(\sigma_3(2^{\lambda_2}) -4\sigma_3(2^{\lambda_2-1})\big) \big(\sigma_3(3^{\lambda_3}) + 9 \sigma_3(3^{\lambda_3-1})\big) \sum_{d\vert m} \mu(d) \left(\frac{3}{d}\right) d \sigma_3(m/d)\\
& = 8 \left(\frac{2^{3\lambda_2+2}+3}{2^3-1}\right) \left(\frac{36 \times 3^{3\lambda_3} -10}{3^3-1}\right)~s_3(m),
\end{split}
\end{equation*}
where $s_3(m)$ is given by \eqref{s3}. The last expression is nothing but Theorem 1.3 of \cite{Ye1}. Now we take ${\bf a} = (1^3,2,3)$. Writing $n$ as above with $\gcd(m,6)=1$, 
we get 
\begin{equation*}
\begin{split}
r_5(1^3, 2, 3;n^2) & = \sum_{d \vert n} \mu(d)d\left(\frac{6}{d}\right) \left\{6 \sigma_3(n/d)-12 \sigma_3(n/2d) -54 
\sigma_3(n/3d)+96 \sigma_3(n/4d)\right.\\
& \left. \qquad \qquad +108  \sigma_3(n/6d)-864 \sigma_3(n/12d) \right\}\\
& = 6 \left(\sigma_3(2^{\lambda_2}) - 2 \sigma_3(2^{\lambda_2-1}) + 16 \sigma_3(2^{\lambda_2-2})\right)
\left(\sigma_3(3^{\lambda_3}) - 9 \sigma_3(3^{\lambda_3-1})\right) \\
& \qquad \qquad \times \sum_{d\vert m} \mu(d) \left(\frac{6}{d}\right) d ~\sigma_3(m/d). 
\end{split}
\end{equation*}
Simplifying the above expression, we get 
\begin{equation*}
r_5(1^3, 2, 3;n^2)  = 12 \left\vert\frac{2^{3\lambda_2+3} -15}{2^3-1}\right\vert \left(\frac{3^{3\lambda_3+2}+4}{3^3-1}\right)
~s_4(m),
\end{equation*}
where $s_4(m)$ is given by \eqref{s4}. The above formula is nothing but Theorem 1.4 of \cite{Ye1}. 

\bigskip

\subsubsection{\bf Data for the remaining cases in Table 1}
In Table 1 there are 54 cases mentioned, out of which we have given explicit formulas for the 18 cases in \corref{quinary}. 
In this section, we obtain the remaining 36 formulas.  In the next table, we give the explicit constants 
$\lambda_{5,1;j}({\bf a})$ that appear in \eqref{eq:main}. These constants are with respect to the corresponding basis elements attached to each vector ${\bf a}$. Using these constants, one gets explicit formulas for the remaining 36 cases. 

The first table gives the linear combination coefficients for the three forms $(1^3,2,4),$ $(1^2,2,4)$, $(1,2,4^3)$. 
We use the basis for $M_4(8)$ given in Table 2, which has dimension 5. The constants are denoted as $\lambda_{4,8;j}$, $1\le j\le 5$.\\

\smallskip

\centerline{\bf Table 4}
\smallskip

\begin{center}
\begin{tabular}{|c|c|c|c|c|c|} 
\hline 
$ {\bf a} $& $\lambda_{4,8;1}$ & $\lambda_{4,8;2}$ & $\lambda_{4,8;3}$ & $\lambda_{4,8;4}$ & $\lambda_{4,8;5}$\\
\hline
$(1^3,2,4)$&$\frac{1}{60}$&$\frac{-1}{60}$&$\frac{1}{30}$&$\frac{-8}{15}$&2\\
\hline 
$(1^2,2,4^2)$&$\frac{1}{120}$&$\frac{-1}{120}$&$\frac{1}{30}$&$\frac{-8}{15}$&2\\
\hline 
$(1,2,4^3)$&$\frac{1}{240}$&$\frac{-1}{240}$&$\frac{1}{30}$&$\frac{-8}{15}$&1\\
\hline 
\end{tabular}
\end{center}

\smallskip

Formulas for $r_5({\bf a};n^2)$ using the above table are given in the following corollary. \\

\begin{cor}\label{cor:4.2}
For a natural number $n= 2^{\lambda_2}m$, where $m$ is an odd positive integer, we have 
\begin{align}
r_5(1^3,2,4;n^2) & = c_{\lambda_2} s_2(m) + 2 C_1(m),\\
r_5(1^2,2,4^2;n^2) & = c'_{\lambda_2} s_2(m) + 2 C_1(m),\\
r_5(1,2,4^3;n^2) & = c''_{\lambda_2} s_2(m) + 2 C_1(m),
\end{align}
where 
$s_2(m)$ is given by \eqref{s2} and $C_1(m)$ is defined as follows:
\begin{equation}
C_1(m) = \tau_{4,8}(m) \prod_{p\ge 3} \left(1 - p \left(\frac{2}{p}\right) \frac{\tau_{4,8}(m/p)}{\tau_{4,8}(m)}\right),
\end{equation}
with $\tau_{4,8}(m)$ being the $m$-th Fourier coefficient of the newform $\Delta_{4,8}(z)$. 
The constants $c_{\lambda_2}$, $c'_{\lambda_2}$ and $c''_{\lambda_2}$ are given in the following table: \\

\begin{center}
\begin{tabular}{|c|c|c|c|}
\hline 
& $\lambda_2=0$ & $\lambda_2=1$ & $\lambda_2\ge 2$\\
\hline
$c_{\lambda_2}$ & $4$ & $32$ &  $12\displaystyle{\left(\frac{9 \times 2^{3\lambda_2 -2} + 10}{2^3-1}\right)}$\\
\hline 
$c'_{\lambda_2}$ & $2$ & $16$ &  $2\displaystyle{\left(\frac{13 \times 2^{3\lambda_2 -1} + 60}{2^3-1}\right)}$\\
\hline 
$c''_{\lambda_2}$ & $1$ & $8$ &  $6\displaystyle{\left(\frac{2^{3\lambda_2} + 20}{2^3-1}\right)}$\\
\hline 
\end{tabular}
\end{center}
\end{cor}
\smallskip

\bigskip

For the remaining 33 cases, we use the basis $M_4(12)$ from Table 2, and in Table 5 below we give the corresponding linear combination coefficients ($\lambda_{4,12;j}$, $1\le j\le 9$) in these cases. \\

\bigskip

\centerline{\bf Table 5}

\bigskip

\begin{center}
\begin{tabular}{|c|c|c|c|c|c|c|c|c|c|}
\hline 
$ {\bf a} $&$\lambda_{4,12;1}$&$\lambda_{4,12;2}$&$\lambda_{4,12;3}$&$\lambda_{4,12;4}$&$\lambda_{4,12;5}$&$\lambda_{4,12;6}$&$\lambda_{4,12;7}$&$\lambda_{4,12;8}$&$\lambda_{4,12;9}$\\
\hline
$(1^3,3^2)$&$\frac{11}{600}$&$\frac{-2}{75}$&$\frac{33}{200}$&0&$\frac{-6}{25}$&0&$\frac{8}{5}$&0&0\\
\hline
$(1,3^2,6^2)$&$\frac{1}{200}$&$\frac{1}{300}$&$\frac{-11}{200}$&0&$\frac{-11}{300}$&0&$\frac{4}{5}$&0&0\\
\hline
$(1^3,6^2)$&$\frac{1}{120}$&$\frac{-1}{60}$&$\frac{3}{40}$&0&$\frac{-3}{20}$&0&4&0&0\\
\hline
$(1,6^4)$&$\frac{1}{600}$&$\frac{1}{150}$&$\frac{-11}{600}$&0&$\frac{-11}{150}$&0&$\frac{8}{5}$&0&0\\
\hline
$(1,2^2,6^2)$&$\frac{1}{200}$&$\frac{-1}{75}$&$\frac{9}{200}$&0&$\frac{-3}{25}$&0 &$\frac{4}{5}$&0&0\\
\hline 
$(2^3,3,6)$&$\frac{1}{300}$&$\frac{1}{75}$&$\frac{-1}{50}$&0 &$\frac{-2}{25}$&0&$\frac{-4}{5}$&0&0\\
\hline 
$(2,3^3,6)$&0&0&$\frac{1}{12}$&0 &$\frac{-1}{6}$&0&0&0&0\\
\hline 
$(2,3,6^3)$&0&0&$\frac{1}{20}$&0 &$\frac{-2}{15}$&0&0&0&0\\
\hline 
$(1^2,2,3,6)$&$\frac{1}{100}$&$\frac{1}{150}$&$\frac{-3}{50}$&0&$\frac{-1}{25}$&0&$\frac{8}{5}$&0&0\\
\hline 
$(2^3,3^2)$&$\frac{1}{120}$&$\frac{-1}{40}$&$\frac{-3}{40}$&$\frac{4}{15}$&$\frac{9}{40}$&$\frac{-12}{5}$&0&0&-2\\
\hline 
$(2,3^4)$&0&0&$\frac{2}{15}$&0 &0&$\frac{-32}{15}$&0&0&0\\
\hline 
%
$(1^3,3,6)$&$\frac{1}{120}$&0&$\frac{1}{8}$&$\frac{-2}{15}$&0&-2&0&0&4\\
\hline 
$(1,3,6^3)$&$\frac{1}{240}$&$\frac{-1}{80}$&$\frac{-17}{240}$&$\frac{2}{15}$&$\frac{17}{80}$&$\frac{-34}{15}$&0&0&1\\
\hline
$(1^2,2,6^2)$&$\frac{1}{120}$&$\frac{-1}{40}$&$\frac{-3}{40}$&$\frac{4}{15}$&$\frac{9}{40}$&$\frac{-12}{5}$&0&0&2\\
\hline 
$(2,3^2,6^2)$&0&0&$\frac{1}{15}$&0&$\frac{1}{15}$&$\frac{-32}{15}$&0&0&0\\
\hline 
$(1,2^2,3,6)$&$\frac{1}{240}$&$\frac{1}{240}$&$\frac{1}{16}$&$\frac{-2}{15}$&$\frac{1}{16}$&-2&0&0&1\\
\hline 
$(1^2,3^3)$&$\frac{1}{100}$&$\frac{1}{25}$&$\frac{-21}{100}$&0&$\frac{-21}{25}$&0&$\frac{8}{5}$&0&0\\
\hline 
$(2^4,3)$&$\frac{1}{150}$&0&$\frac{3}{50}$&$\frac{-8}{75}$&0&$\frac{-24}{25}$&$\frac{-8}{5}$&$\frac{-32}{5}$&0\\
\hline 
$(2^2,3^3)$&$\frac{1}{150}$&$\frac{-1}{100}$&$\frac{-7}{50}$&$\frac{4}{75}$&$\frac{21}{100}$&$\frac{-28}{25}$&$\frac{-4}{15}$&$\frac{-128}{15}$&$\frac{-4}{3}$\\
\hline 
$(1^2, 2^2,3)$&$\frac{1}{75}$&$\frac{-1}{150}$&$\frac{3}{25}$&$\frac{-8}{75}$&$\frac{-3}{50}$&$\frac{-24}{25}$&$\frac{4}{5}$&$\frac{32}{5}$&0\\
\hline 
$(1^2,3,6^2)$&$\frac{1}{150}$&$\frac{-1}{100}$&$\frac{-7}{50}$&$\frac{4}{75}$&$\frac{21}{100}$&$\frac{-28}{25}$&$\frac{16}{15}$&$\frac{32}{15}$&$\frac{4}{3}$\\
\hline 
$(1^3,2,6)$&$\frac{1}{60}$&$\frac{-1}{40}$&$\frac{-3}{20}$&$\frac{2}{15}$&$\frac{9}{40}$&$\frac{-6}{5}$&0&0&2\\
\hline 
$(1,2,6^3)$&$\frac{1}{600}$&0&$\frac{13}{200}$&$\frac{-2}{75}$&0&$\frac{-26}{25}$&$\frac{4}{15}$&$\frac{-16}{15}$&$\frac{4}{3}$\\
\hline 
$(2^2,3,6^2)$&$\frac{1}{300}$&$\frac{-1}{150}$&$\frac{-7}{100}$&$\frac{4}{75}$&$\frac{7}{50}$&$\frac{-28}{25}$&$\frac{-4}{5}$&$\frac{-16}{5}$&0\\
\hline 
$(1,2,3^2,6)$&$\frac{1}{300}$&$\frac{-1}{600}$&$\frac{13}{100}$&$\frac{-2}{75}$&$\frac{-13}{200}$&$\frac{-26}{25}$&$\frac{8}{15}$&$\frac{64}{15}$&0\\
\hline
$(1,2^3,3)$&$\frac{1}{80}$&$\frac{-3}{80}$&$\frac{-9}{80}$&$\frac{2}{5}$&$\frac{27}{80}$&$\frac{-18}{5}$&0&0&-1\\
\hline 
$(1,2,3^3)$&$\frac{1}{200}$&0&$\frac{39}{200}$&$\frac{-2}{25}$&0&$\frac{-78}{25}$&$\frac{-8}{15}$&$\frac{32}{15}$&$\frac{4}{3}$\\
\hline 
$(1^4,6)$&$\frac{1}{50}$&0&$\frac{9}{50}$&$\frac{-8}{25}$&0&$\frac{-72}{25}$&$\frac{16}{5}$&$\frac{-64}{5}$&0\\
\hline 
$(1^2,6^3)$&$\frac{1}{200}$&$\frac{-3}{200}$&$\frac{-21}{200}$&$\frac{4}{25}$&$\frac{63}{200}$&$\frac{-84}{25}$&$\frac{32}{15}$&$\frac{16}{15}$&$\frac{2}{3}$\\
\hline 
$(1^2,3^2,6)$&$\frac{1}{100}$&$\frac{-1}{50}$&$\frac{-21}{100}$&$\frac{4}{25}$&$\frac{21}{50}$&$\frac{-84}{25}$&$\frac{8}{5}$&$\frac{32}{5}$&0\\
\hline 
$(1^2,2^2,6)$&$\frac{1}{100}$&$\frac{1}{100}$&$\frac{9}{100}$&$\frac{-8}{25}$&$\frac{9}{100}$&$\frac{-72}{25}$&$\frac{8}{5}$&$\frac{16}{5}$&0\\
\hline 
$(2^2,3^2,6)$&$\frac{1}{200}$&$\frac{-3}{200}$&$\frac{-21}{200}$&$\frac{4}{25}$&$\frac{63}{200}$&$\frac{-84}{25}$&$\frac{-8}{15}$&$\frac{-64}{15}$&$\frac{-2}{3}$\\
\hline
$(1,2,3,6^2)$&$\frac{1}{400}$&$\frac{1}{400}$&$\frac{39}{400}$&$\frac{-2}{25}$&$\frac{39}{400}$&$\frac{-78}{25}$&$\frac{16}{15}$&$\frac{32}{15}$&$\frac{1}{3}$\\
 \hline
\end{tabular}
\end{center}

\bigskip

\noindent 
Out of these 33 cases, 11 can be given in simpler form. These 11 cases are: $(1^3,3^2), (1,3^2,6^2),$ $(1^3,6^2), (1, 6^4), (1, 2^2,6^2), (2^3, 3,6), (2, 3^3,6), (2,3,6^3),$ $(1^2,2,3,6), (2,3^4),$ $(2,3^2,6^2)$.  We list the formulas corresponding to these cases in the following corollary.

\bigskip

\begin{cor}\label{cor:4.3}
Let $n$ be a natural number and write it as $n = 2^{\lambda_2} 3^{\lambda_3} m$, where $\gcd(m,6) =1$. Then we have the following formulas (using data from Table 5):

\begin{align}
r_5({\bf a}; n^2) & = c_{{\bf a},\lambda_2} d_{\lambda_3} s_1(m) + e_{\bf a} \tau_{4,6}(2^{\lambda_2} 3^{\lambda_3}) C_2(m),\\
r_5(2,3^4;n^2) & = 96 \left(\frac{2^{3\lambda_2+1}+5}{2^3-1}\right) \sigma_3(3^{\lambda_3-1}) s_2(m),\\
r_5(2,3^2,6^2;n^2) & =  80 \left(\frac{2^{3\lambda_2}+6}{2^3-1}\right) \sigma_3(3^{\lambda_3-1}) s_2(m),
\end{align}
where $s_1(m)$, $s_2(m)$ are defined by \eqref{s1}, \eqref{s2} and for $m\ge 1$, $\gcd(m,6)=1$, $C_2(m)$ is defined by 

\begin{equation}
C_2(m) = \tau_{4,6}(m) \prod_{p\ge 5}\left(1 - p\frac{\tau_{4,6}(m/p)}{\tau_{4,6}(m)}\right).
\end{equation}
The constants $c_{{\bf a},\lambda_2}$, $d_{\lambda_3}$ and $e_{\bf a}$ for the 9 vectors ${\bf a}$ are given in the table 
below. \\

\begin{center}
\begin{tabular}{|c|c|c|c|}
\hline 
{\bf a} & $c_{{\bf a},\lambda_2}$ & $d_{\lambda_3}$ & $e_{\bf a}$\\
\hline 
$(1^3,3^2)$ & $\displaystyle{\frac{4}{5} \left( \frac{9 \times 2^{3\lambda_2+3} + 5}{2^3-1}\right)}$ & 
$ \displaystyle{\left(\frac{2 \times 3^{3\lambda_3+2} -5 }{3^3-1}\right)}$ & $\dfrac{8}{5}$\\
\hline 
$(1,3^2,6^2)$ & $\displaystyle{\frac{2}{5} \left(\frac{13 \times 2^{3\lambda_2+1} - 5}{2^3-1}\right)}$ & 
$ \displaystyle{\left(\frac{16 \times 3^{3\lambda_3} +10 }{3^3-1}\right)}$ & $\dfrac{4}{5}$\\
\hline 
$(1^3,6^2)$ &$\displaystyle{2 \left(\frac{3 \times 2^{3\lambda_2+1} +1}{2^3-1}\right)}$ & 
$\displaystyle{\left(\frac{4 \times 3^{3\lambda_3+2} -10 }{3^3-1}\right)}$  & $4$\\
\hline
$(1,6^4)$& $\displaystyle{\frac{2}{5} \left(\frac{3 \times 2^{3\lambda_2+2} -5}{2^3-1}\right)}$ & 
$\displaystyle{\left(\frac{16 \times 3^{3\lambda_3} +10 }{3^3-1}\right)}$ & $\dfrac{8}{5}$\\
\hline 
$(1,2^2,6^2)$ & $\displaystyle{\frac{2}{5} \left(\frac{2^{3\lambda_2+4} + 5}{2^3-1}\right)}$ & 
$ \displaystyle{\left(\frac{4 \times 3^{3\lambda_3+2} -10 }{3^3-1}\right)}$ & $\dfrac{4}{5}$\\
\hline 
$(2^3,3,6)$ & $\displaystyle{\frac{4}{5} \left(\frac{3 \times 2^{3\lambda_2+2} - 5}{2^3-1}\right)}$ & 
$ \displaystyle{\left(\frac{7\times 3^{3\lambda_3+1} +5 }{3^3-1}\right)}$ & $\dfrac{-4}{5}$\\
\hline 
$(2,3^3,6)$ & $\displaystyle{20\left(\frac{3\times 2^{3\lambda_2+1}+1}{2^3-1}\right)}$ & $\sigma_3(3^{\lambda_3-1})$ &$0$\\
\hline 
$(2,3,6^3)$ & $\displaystyle{4\left(\frac{2^{3\lambda_2+4}+5}{2^3-1}\right)}$ & $\sigma_3(3^{\lambda_3-1})$ &$0$\\
\hline  
$(1^2,2,3,6)$ & $\displaystyle{\frac{4}{5} \left(\frac{13\times 2^{3\lambda_2+1} - 5}{2^3-1}\right)}$ & 
$ \displaystyle{\left(\frac{11 \times 3^{3\lambda_3+1} -7 }{3^3-1}\right)}$ & $\dfrac{8}{5}$\\
\hline 
\end{tabular}
\end{center}
\end{cor}

\bigskip

\noindent {\bf Note}: In the above corollary, $\lambda_3\ge 1$ when ${\bf a} = (2,3^3,6)$ and $(2,3,6^3)$. 
\bigskip

 \subsection{Septenary forms}
 
In this section we discuss the septenary case $\ell =7$. We find formulas when the coefficients $a_i$ belong to the set 
$\{1,2,4\}$. There are 27 forms corresponding to the choice of $a_i$'s in this set. Let $i_1 + i_2+i_3 =7$. i.e., the coefficient 1 appears $i_1$ times, 2 appears $i_2$ times and $4$ appears $i_3$ times. The theta series associated to these 27 forms 
are given by  $\Theta^{i_1}(z) \Theta^{i_2}(2z) \Theta^{i_3}(4z)$. They belong to the one of the spaces $M_{7/2}(8)$, $M_{7/2}(8,\chi_2)$,  $M_{7/2}(16)$, $M_{7/2}(16,\chi_2)$ according as $i_3=0, i_2$ even; $i_3=0, i_2$ odd; $i_3\not=0, i_2$ even; and $i_3\not=0, i_2$ odd respectively. Since $(\ell-1)/2$ is odd in this case, we take $D =-4$ and apply the Shimura map ${\mathcal S}_D$ on these theta products and by using \thmref{thm:main} the images of these functions are in $M_6(4)$ or $M_6(8)$ according as $i_3=0$ or not.  We now use the following basis for these spaces: 
 
\bigskip
 
 \centerline{List of basis for $M_6(4)$ and $M_6(8)$}
 
\begin{center}
\smallskip
\smallskip
\begin{tabular}{|c|c|c|}
\hline 
Space & Basis elements & Dimension\\
\hline 
$M_6(4)$ & $E_6(z),E_6(2z),E_6(4z),\Delta_{6,4}(z)$ & 4\\
\hline 
$M_6(8)$ & $E_6(z),E_6(2z),E_6(4z),E_6(8z), \Delta_{6,4}(z),  \Delta_{6,4}(2z), \Delta_{6,8}(z)$&  7\\
\hline 
\end{tabular}
\end{center}

\smallskip
\bigskip
The newforms $\Delta_{k,N}(z)$ that appear in the above table are defined below. 

\begin{equation*}
\begin{split}
\Delta_{6,4}(z) & = \eta^{12}(2z) = \sum_{n\ge 1} \tau_{6,4}(n) q^n,\\
\Delta_{6,8}(z) & = \frac{4}{45} \left[ \frac{1}{5376} E_6(z) + \frac{5}{1792} E_6(2z) + \frac{5}{112} E_6(4z)+ \frac{20}{21} E_6(8z) - \frac{405}{32}\Delta_{6,4}(z) \right.\\
& \qquad \quad \left.-\frac{135}{2}\Delta_{6,4}(2z) - E_2(z)E_4(8z)+ 3 DE_4(8z)\right] = \sum_{n\ge 1} \tau_{6,8}(n) q^n.
\end{split}
\end{equation*}
We  have another expression for $\Delta_{6,8}(z)$, which we give below.
\begin{equation*}
\begin{split}
\Delta_{6,8}(z) &= \frac{-1}{90} \left[ \frac{5}{336} E_6(z) + \frac{5}{112} E_6(2z) + \frac{5}{28} E_6(4z)+ \frac{16}{21} E_6(8z) + \frac{135}{2}\Delta_{6,4}(z) \right.\\
& \qquad \quad \left.+810\Delta_{6,4}(2z) - E_2(8z)E_4(z)+ \frac{3}{8} DE_4(z)\right].
\end{split}
\end{equation*}
In the above, $Df$ denotes the derivative of the modular form $f$ with respect to $z$, which is given by $\frac{1}{2\pi i} f'(z)$. 
It is known that $Df(z) - \frac{k}{12} E_2(z) f(z)$ is a modular form of weight $k+2$, where $f$ is a modular form of weight $k$ and $E_2(z) = 1 -24 \sum_{n\ge 1} \sigma(n) q^n$ is the normalised Eisenstein series of weight $2$ (which is a quasimodular form), where $\sigma(n)$ is the sum of the positive divisors of $n$. 

\smallskip

We use the basis elements given in the above table to obtain the following explicit expressions for ${\mathcal S}_{-4}(\Theta_{\bf a}(z))$, where ${\bf a} = (1^{i_1}, 2^{i_2}, 4^{i_3})$, $i_1 + i_2+i_3 =7$, with either $i_3=0$ or $i_3\not= 0$, $i_2$ is even. There are 18 vectors ${\bf a}$ with these conditions. \\

\smallskip

\centerline{\bf Table 6}

\smallskip

\begin{center}
{\small 
\begin{tabular}{|c|c|}
\hline 
${\bf a}$ & ${\mathcal S}_{-4}(\Theta_{\bf a})(z)$\\
\hline 
$(1^6,2)$ & $\frac{-31}{42} E_6(z) - \frac{16}{21} E_6(2z)$\\
\hline 
$(1^5,2^2)$ &$\frac{-127}{252} E_6(z) + \frac{16}{63} E_6(2z)$\\
\hline 
$(1^4,2^3)$ &$\frac{-5}{14} E_6(z) - \frac{8}{7} E_6(2z)$\\
\hline 
$(1^3,2^4),$&\\
$(1^4,2^2,4)$ &$ \frac{-1}{4} E_6(z)$ \\
\hline 
$(1^2,2^5)$ &$\frac{-1}{6} E_6(z) - \frac{4}{3} E_6(2z)$\\
\hline 
$(1^6,4)$ &$ \frac{-127}{252} E_6(z) + \frac{16}{63} E_6(2z)$\\
\hline 
\end{tabular}
\quad 
\begin{tabular}{|c|c|}
\hline 
${\bf a}$ & ${\mathcal S}_{-4}(\Theta_{\bf a})(z)$\\
\hline 
$(1^5,4^2)$ &$\frac{-47}{252} E_6(z) - \frac{4}{63} E_6(2z)$\\
\hline 
$(1^3,4^4), (1^2,4^5), $ &\\
$(1,4^6), (1,2^2,4^4)$ & $ \frac{-1}{36} E_6(z) - \frac{2}{9} E_6(2z)$\\
\hline
$(1^3,2^2,4^2),$ & \\
$(1^2,2^4,4), (1,2^6)$ & $\frac{-31}{252} E_6(z) - \frac{8}{63} E_6(2z)$\\
\hline 
$(1^2,2^2,4^3),$ &\\
$(1,2^4,4^2), (1^4,4^3) $ & $\frac{-5}{84} E_6(z) - \frac{4}{21} E_6(2z)$\\
\hline 
\end{tabular}
}
\end{center}

\smallskip

\smallskip

\smallskip
Using the above data, we deduce some of the results and identities proved in \cite{Ye2}. Put $n= 2^km$, $k\ge 0, m\ge 1$ are integers, $m$ odd. 
Comparing the Fourier coefficients and simplifying (using the multiplicative property of the divisor function $\sigma_5(n)$)
as done earlier, we get the following formulas as corollary to our main theorem. When $n$ is even, the formulas given below 
for $r_7({\bf a}; 4n^2)$ (derived from our method) are the same 18 formulas obtained in Theorems 2.1 to 2.3 in \cite{Ye2}.

\begin{cor}\label{3}
For an integer $k\ge 1$, we have the following identities:
\begin{align}
r_7(1^6,2;2^{2k}m^2) & = 12 \left(\frac{2^{5k+5} -63}{2^5-1}\right) s_5(m),\label{22}\\
r_7(1^5,2^2;2^{2k}m^2) & = ~r_7(1^6,4;2^{2k}m^2) ~=~ \left(\frac{250 \times 2^{5k} -126}{2^5-1}\right) s_6(m),\label{23}\\
r_7(1^4,2^3;2^{2k}m^2) & = \left(\frac{198 \times 2^{5k} -756}{2^5-1}\right) s_5(m),\label{24}
\end{align}
\begin{align}
r_7(1^3,2^4; 2^{2k}m^2) & = r_7(1^4,2^2,4;2^{2k}m^2) ~=~126 \left(\frac{2^{5k}-1}{2^5-1}\right) s_6(m),\label{25}\\
r_7(1^2,2^5;2^{2k}m^2) & =  \left(\frac{105 \times 2^{5k} -756}{2^5-1}\right) s_5(m),\label{26}\\
r_7(1^5,4^2;2^{2k}m^2) & = \left(\frac{95\times 2^{5k} -126}{2^5-1}\right) s_6(m),\label{27}\\
r_7(1^3,4^4; 2^{2k}m^2) & = r_7(1^2,4^5; 2^{2k}m^2) ~=~ r_7(1,4^6; 2^{2k}m^2) \nonumber \\
& =~ r_7(1,2^2,4^4;2^{2k}m^2)~ = ~\left(\frac{35\times 2^{5k-1} -126}{2^5-1}\right) s_6(m),\label{28}\\
r_7(1^3,2^2,4^2;2^{2k}m^2) & = r_7(1^2,2^4,4;2^{2k}m^2) ~=~ r_7(1,2^6;2^{2k}m^2)  \nonumber \\
& = \left(\frac{2^{5k+6} -126}{2^5-1}\right) s_6(m),\label{29}\\
r_7(1^2,2^2,4^3;2^{2k}m^2) &= r_7(1,2^4,4^2;2^{2k}m^2)  = r_7(1^4,4^3;2^{2k}m^2)  \nonumber \\
&=  ~\left(\frac{33\times 2^{5k} -126}{2^5-1}\right) s_6(m),\label{30}
\end{align}
where $s_5(m)$ and $s_6(m)$ are given by \eqref{sm} and \eqref{tm}.
\end{cor}

\smallskip

\begin{rmk}
Formulas \eqref{22}, \eqref{24}, \eqref{26} given above are exactly the same as in Lemma 3.1 of \cite{Ye2}, except for the fact that $r_7(1^6,2;m^2) = 12 s_5(m)$ (see \rmkref{rmk:conj}).  Formulas \eqref{29} (${\bf a} =(1,2^6)$), \eqref{23} (${\bf a} 
= (1^6,4)$), \eqref{28} (${\bf a} = (1,4^6)$) are the same as in  Lemmas 4.1, 4.2 and Lemma 5.1 of \cite{Ye2}, except for the fact that $r_7(m^2) = 14 s_6(m)$. We also note that the same formula for $r_7(n^2)$ can be derived using the Shimura correspondence and this fact was observed in \cite[\S 3, p.372]{gun}. 
We also note that the above corollary proves all the results stated  in Theorems 2.1, 2.2, 2.3 of \cite{Ye2}, when $n$ is an even integer.   
\end{rmk}

\bigskip

\subsubsection{\bf Remaining 9 cases} ~We now consider the remaining 9 formulas, where $i_3\not= 0$ and $i_2$ is odd and obtain new formulas.  We again take $D =-4$ and use  \thmref{thm:main} to obtain the following explicit images under the Shimura lifting:

\begin{equation*}
\begin{split}
{\mathcal S}_{-4}\left(\Theta(z) \Theta(2z)\Theta^5(4z)\right) & = \frac{-1}{42} E_6(z) + \frac{1}{21} E_6(2z) - \frac{32}{21} E_6(4z),\\
{\mathcal S}_{-4}\left(\Theta(z) \Theta^3(2z)\Theta^3(4z)\right) & = {\mathcal S}_{-4}\left(\Theta^2(z) \Theta(2z)\Theta^4(4z)\right) \\
& = \frac{-1}{21} E_6(z) + \frac{1}{14} E_6(2z) - \frac{32}{21} E_6(4z) - 4 \Delta_{6,4}(z),\\
{\mathcal S}_{-4}\left(\Theta(z) \Theta^5(2z)\Theta(4z)\right) & = {\mathcal S}_{-4}\left(\Theta^2(z) \Theta^3(2z)\Theta^2(4z)\right) \\
& = \frac{-2}{21} E_6(z) + \frac{5}{42} E_6(2z) - \frac{32}{21} E_6(4z) - 4 \Delta_{6,4}(z),
\end{split}
\end{equation*}
\begin{equation*}
\begin{split}
{\mathcal S}_{-4}\left(\Theta^3(z) \Theta(2z)\Theta^3(4z)\right)
& = \frac{-2}{21} E_6(z) + \frac{5}{42} E_6(2z) - \frac{32}{21} E_6(4z) - 12 \Delta_{6,4}(z),\\
{\mathcal S}_{-4}\left(\Theta^3(z) \Theta^3(2z)\Theta(4z)\right)
& = \frac{-4}{21} E_6(z) + \frac{3}{14} E_6(2z) - \frac{32}{21} E_6(4z) - 4 \Delta_{6,4}(z),\\
{\mathcal S}_{-4}\left(\Theta^4(z) \Theta(2z)\Theta^2(4z)\right)
& = \frac{-4}{21} E_6(z) + \frac{3}{14} E_6(2z) - \frac{32}{21} E_6(4z) - 20 \Delta_{6,4}(z),\\
{\mathcal S}_{-4}\left(\Theta^5(z) \Theta(2z)\Theta(4z)\right)
& = \frac{-8}{21} E_6(z) + \frac{17}{42} E_6(2z) - \frac{32}{21} E_6(4z) - 20 \Delta_{6,4}(z).\\
\end{split}
\end{equation*}

\smallskip

\bigskip

\noindent 
From the above, we deduce the following 9 formulas which are new (and not obtained in the work of Cooper-Lam-Ye \cite{Ye2}). 
Write $n = 2^{\lambda_2}m$, $\lambda_2\ge 1$, $2\not\!\vert ~m$. Let $m= \displaystyle{\prod_{p {\rm ~odd~}}} p^{\lambda_p}$, $\lambda_p \ge 0$ and let $s_5(m)$ be as defined in \eqref{sm}. It is a fact that $\tau_{6,4}(n) = 0$, when $n$ is an even positive integer and so, we define the function $C_3(m)$ as follows: 
\begin{equation}\label{c1}
C_3(m) = \tau_{6,4}(m) \prod_{p\ge 3} \left(1 - p^2 \left(\frac{-2}{p}\right) \frac{\tau_{6,4}(m/p)}{\tau_{6,4}(m)}\right).
\end{equation}

Using the multiplicative property of $\sigma_5(n)$ and $\tau_{6,4}(n)$, we get the following new formulas (as a consequence of the explicit Shimura images stated above):

\begin{cor}\label{even}
For a natural number $n= 2^{\lambda_2} m$, $\lambda_2\ge 1$, $m\ge 1$ odd, we have 
\begin{equation}
r_7(1,2,4^5;n^2) = 12 \left\vert\displaystyle{\frac{2^{5\lambda_2}- 63}{2^5-1}}\right\vert s_5(m),\hskip 3.25cm 
\end{equation}

\begin{equation}
r_7(1,2^3,4^3;n^2) = r_7(1^2,2,4^4;n^2) = \begin{cases} 24 s_5(m) - 4 C_3(m) &  {\rm ~if~} \lambda_2 =1, \\
\displaystyle{756 \left(\frac{2^{5\lambda_2-5} -1}{2^5-1}\right)} s_5(m) &  {\rm ~if~} \lambda_2 \ge 2, \\
\end{cases}
\end{equation}

\begin{equation}
r_7(1,2^5,4; n^2) = r_7(1^2,2^3,4^2; n^2) =\begin{cases} 48 s_5(m) - 4 C_3(m) &  {\rm ~if~} \lambda_2 =1, \\
\displaystyle{\frac{375 \times 2^{5\lambda_2-3} -756}{2^5-1}} s_5(m) &  {\rm ~if~} \lambda_2 \ge 2, \\
\end{cases}
\end{equation}

\begin{equation}
r_7(1^3,2,4^3; n^2) = \begin{cases} 48 s_5(m) - 12 C_3(m) &  {\rm ~if~} \lambda_2 =1, \\
\displaystyle{\frac{375 \times 2^{5\lambda_2-3} -756}{2^5-1}} s_5(m) &  {\rm ~if~} \lambda_2 \ge 2, \\
\end{cases}
\end{equation}

\begin{equation}
r_7(1^3,2^3,4; n^2) = \begin{cases} 96 s_5(m) - 4 C_3(m) &  {\rm ~if~} \lambda_2 =1, \\
\displaystyle{\frac{747 \times 2^{5\lambda_2-3} -756}{2^5-1}} s_5(m) &  {\rm ~if~} \lambda_2 \ge 2,\\
\end{cases}
\end{equation}

\begin{equation}
r_7(1^4,2,4^2; n^2) = \begin{cases} 96 s_5(m) - 20 C_3(m) &  {\rm ~if~} \lambda_2 =1, \\
\displaystyle{\frac{747 \times 2^{5\lambda_2-3} -756}{2^5-1}} s_5(m) &  {\rm ~if~} \lambda_2 \ge 2, \\
\end{cases}
\end{equation}

\begin{equation}
r_7(1^5,2,4; n^2) = \begin{cases} 192 s_5(m) - 20 C_3(m) &  {\rm ~if~} \lambda_2 =1, \\
\displaystyle{\frac{2982 \times 2^{5\lambda_2-4} -756}{2^5-1}} s_5(m) &  {\rm ~if~} \lambda_2 \ge 2. \\
\end{cases}
\end{equation}
\end{cor}

\smallskip

\begin{rmk}\label{new1}
As mentioned earlier, we could use only the Shimura map ${\mathcal S}_{-4}$ in the case when 
$\ell =7$ (more generally $\ell \equiv 3\pmod{4}$) in order to  apply \thmref{thm:main}. So, we could get formulas only 
for $4n^2$ instead of $n^2$. In this connection, we refer to the concluding remarks in \cite{Ye2}, where they also mention that these formulas (for $4n^2$) can be obtained using methods of \cite{cooper-lam}. However, explicit formulas for these cases have not been given in their work. The above corollary gives explicit formulas in this case and they are new. In the next section, we determine these formulas for $n^2$ for any natural number $n\ge 1$ (under an assumption). 
\end{rmk}

\begin{rmk}
Using the formulas given in the above corollary, we obtain several relations connecting $r_7({\bf a};n^2)$ and $C_3(m)$. Some of them are listed in the following Proposition. It is possible that these relations can be obtained using different methods. 
\end{rmk}

\smallskip

\begin{prop}
Let $n\ge 1$ be an integer and $m$ be the odd part of $n$. We have the following relations, which arise using the identities listed in \corref{even}. 
When $n\equiv 0\pmod{4}$, we have 
\begin{align}
r_7(1,2^5,4;n^2) & = ~r_7(1^2,2^3,4^2; n^2) ~=~ r_7(1^3,2,4^3;n^2),\\
r_7(1^3,2^3,4;n^2) & = ~r_7(1^4,2,4^2;n^2).
\end{align}
When $n\equiv 2\pmod{4}$, we have 
\begin{align}
4 C_3(m) & = r_7(1,2^5,4;n^2)  - 2~ r_7(1,2^3,4^3; n^2)  = r_7(1^2,2^3,4^2;n^2)  - 2~ r_7(1,2^3,4^3; n^2) 
\nonumber \\
& = r_7(1,2^5,4;n^2)  - 2~ r_7(1^2,2,4^4; n^2)  =  r_7(1^2,2^3,4^2;n^2)  - 2~ r_7(1^2,2,4^4; n^2),\\
8 C_3(m) & = ~ r_7(1,2^5,4;n^2)  - r_7(1^3,2,4^3; n^2)  = ~ r_7(1^2,2^3,4^2;n^2)  - 2 r_7(1^3,2,4^3; n^2),\\
r_7(1,2^5,4;n^2) & =   4 \big(r_7(1^2,2,4^4; n^2) - r_7(1^3,2,4^3; n^2)\big),
\end{align}
where $C_3(m)$ is given by \eqref{c1}.
\end{prop}

\bigskip

\subsubsection{\bf Conjectural formulas} ~ 
In this section, we assume that the Shimura map ${\mathcal S}_{-1}$ maps the space $M_{k+1/2}(4N, \chi)$ into the space $M_{2k}(2N)$, where $\chi$ is a quadratic Dirichlet character modulo $4N$. All the results presented in this section are based on this assumption. Therefore, it implies that  ${\mathcal S}_{-1}(\Theta^{i_1}(z) \Theta^{i_2}(2z) \Theta^{i_3}(4z)) \in M_6(8)$, 
where $(1^{i_1},2^{i_2},4^{i_3}) \in \{(1,2,4^5), (1,2^3,4^3), (1^2,2,4^4), (1,2^5,4), (1^2,2^3,4^2), (1^3,2,4^3), (1^3,2^3,4),
(1^4,2,4^2), (1^5,2,4)\}$ and the explicit images of these functions are given below.
 
 \smallskip

\begin{equation*}
\begin{split}
{\mathcal S}_{-1}\left(\Theta(z) \Theta(2z)\Theta^5(4z)\right) & = \frac{-1}{1344} E_6(z) + \frac{1}{1344} E_6(2z) + \frac{1}{42} E_6(4z) -\frac{32}{21} E_6(8z) \\
&\qquad + \frac{5}{8} \Delta_{6,4}(z) + \Delta_{6,8}(z),\\
{\mathcal S}_{-1}\left(\Theta(z) \Theta^3(2z)\Theta^3(4z)\right) & = 
\frac{-1}{672} E_6(z) + \frac{1}{672} E_6(2z) + \frac{1}{42} E_6(4z) -\frac{32}{21} E_6(8z) \\
& \qquad + \frac{1}{4} \Delta_{6,4}(z)  - 4 \Delta_{6,4}(2z) + \Delta_{6,8}(z),
\end{split}
\end{equation*}
\begin{equation*}
\begin{split}
{\mathcal S}_{-1}\left(\Theta^2(z) \Theta(2z)\Theta^4(4z)\right) & = 
\frac{-1}{672} E_6(z) + \frac{1}{672} E_6(2z) + \frac{1}{42} E_6(4z) -\frac{32}{21} E_6(8z) \\
& \qquad + \frac{5}{4} \Delta_{6,4}(z)  - 4 \Delta_{6,4}(2z) + 2 \Delta_{6,8}(z),\\
{\mathcal S}_{-1}\left(\Theta(z) \Theta^5(2z)\Theta(4z)\right) & = 
\frac{-1}{336} E_6(z) + \frac{1}{336} E_6(2z) + \frac{1}{42} E_6(4z) -\frac{32}{21} E_6(8z) \\
& \qquad + \frac{1}{2} \Delta_{6,4}(z)  - 4 \Delta_{6,4}(2z),\\
{\mathcal S}_{-1}\left(\Theta^2(z) \Theta^3(2z)\Theta^2(4z)\right) &= 
\frac{-1}{336} E_6(z) + \frac{1}{336} E_6(2z) + \frac{1}{42} E_6(4z) -\frac{32}{21} E_6(8z) \\
& \qquad + \frac{1}{2} \Delta_{6,4}(z) - 4 \Delta_{6,4}(2z) + 2 \Delta_{6,8}(z),\\
{\mathcal S}_{-1}\left(\Theta^3(z) \Theta(2z)\Theta^3(4z)\right) &=
\frac{-1}{336} E_6(z) + \frac{1}{336} E_6(2z) + \frac{1}{42} E_6(4z) -\frac{32}{21} E_6(8z) \\
& \qquad - \frac{3}{2} \Delta_{6,4}(z) - 12 \Delta_{6,4}(2z) + 3 \Delta_{6,8}(z),\\
{\mathcal S}_{-1}\left(\Theta^3(z) \Theta^3(2z)\Theta(4z)\right) & = 
\frac{-1}{168} E_6(z) + \frac{1}{168} E_6(2z) + \frac{1}{42} E_6(4z) -\frac{32}{21} E_6(8z) \\
& \quad +  \Delta_{6,4}(z)  - 4 \Delta_{6,4}(2z) + 2 \Delta_{6,8}(z),\\
{\mathcal S}_{-1}\left(\Theta^4(z) \Theta(2z)\Theta^2(4z)\right) & = 
\frac{-1}{168} E_6(z) + \frac{1}{168} E_6(2z) + \frac{1}{42} E_6(4z) -\frac{32}{21} E_6(8z)\\
& \qquad  +  \Delta_{6,4}(z)  - 20 \Delta_{6,4}(2z) + 4 \Delta_{6,8}(z),\\
{\mathcal S}_{-1}\left(\Theta^5(z) \Theta(2z)\Theta(4z)\right) &=
\frac{-1}{84} E_6(z) + \frac{1}{84} E_6(2z) + \frac{1}{42} E_6(4z) -\frac{32}{21} E_6(8z)  \\
& \quad - 20 \Delta_{6,4}(2z) 
+ 4 \Delta_{6,8}(z).\\
\end{split}
\end{equation*}

\smallskip

As done before, we can simplify the above formulas. We write $n = 2^{\lambda_2}m$, with $\lambda_2\ge 0$, $2\not\vert m$ and $m$ is expressed as before. Let $s_5(m)$ be as defined in \eqref{sm}. In this case $\tau_{6,8}(n)$ appear and they are zero when $n$ is even. For $m\ge 1$ odd, let $C_3(m)$ be as defined by \eqref{c1} and we let 
\begin{equation}
C_4(m) = \tau_{6,8}(m) \prod_{p\ge 3} \left(1 - p^2 \left(\frac{-2}{p}\right) \frac{\tau_{6,8}(m/p)}{\tau_{6,8}(m)}\right).
\end{equation}

\noindent 
Using the above explicit expressions given in terms of the basis elements of $M_6(8)$ and using the multiplicative properties of $\sigma_5(n)$, $\tau_{6,4}(n)$, $\tau_{6,8}(n)$, the following (conjectural) formulas are obtained. 

\smallskip
\noindent {\bf Conjectural formulas}. Let $m\ge 1$ be an odd integer. Then, we have 

\begin{eqnarray}
r_7(1,2,4^5; m^2) &=& \frac{3}{8} s_5(m) + \frac{5}{8} C_3(m) +  C_4(m), \\
r_7(1,2^3,4^3; m^2) &=&  \frac{3}{4} s_5(m) + \frac{1}{4} C_3(m) +  C_4(m), \\
r_7(1^2,2,4^4; m^2) &=& \frac{3}{4} s_5(m) + \frac{5}{4} C_3(m) + 2 C_4(m),\\
r_7(1,2^5,4; m^2) &=& \frac{3}{2} s_5(m) + \frac{1}{2} C_3(m),\\ 
r_7(1^2,2^3,4^2; m^2) &=& \frac{3}{2} s_5(m) + \frac{1}{2} C_3(m) + 2 C_4(m),\\
r_7(1^3,2,4^3; m^2) &=&  \frac{3}{2} s_5(m) + \frac{3}{2} C_3(m) + 3 C_4(m),\\
r_7(1^3,2^3,4; m^2) &=& 3 s_5(m) + C_3(m) + 2 C_4(m),\\
r_7(1^4,2,4^2; m^2) &=&  3 s_5(m) + C_3(m) + 4 C_4(m),\\
r_7(1^5,2,4; m^2) &=& 6 s_5(m) + 4 C_4(m).
\end{eqnarray}

\smallskip

\begin{rmk}
When $n$ is an even positive integer, the formulas for $r_7({\bf a};n^2)$ obtained using our assumption (i.e., using the Shimura map ${\mathcal S}_{-1}$) coincide with the formulas proved in \corref{even} (which was obtained using the Shimura map ${\mathcal S}_{-4}$). 
\end{rmk}

\bigskip
\begin{rmk}\label{rmk:conj}
Under the assumption made at the beginning of this section, we also have ${\mathcal S}_{-1}(\Theta^6(z) \Theta(2z)) \in M_6(4)$ and writing it in terms of the basis elements mentioned in \S4.2, we get 
\begin{equation*}
{\mathcal S}_{-1}(\Theta^6(z) \Theta(2z)) = \frac{-1}{42} E_6(z) + \frac{1}{21} E_6(2z) - \frac{32}{21} E_6(4z).
\end{equation*}
Like we did earlier, writing $n = 2^\lambda m$, $\lambda\ge 0$ and $m\ge 1$ odd, the above identity can be simplified 
(after comparing the $n$-th Fourier coefficients and taking M\"obius inversion) and we get the following formula (which is 
exactly the first formula in Theorem 2.1 of \cite{Ye2}): 
\begin{equation}
r_7(1^6,2; n^2) = 12 \left|\frac{2^{5\lambda +5}-63}{2^5-1}\right| s_5(m),
\end{equation}
where $s_5(m)$ is given by \eqref{sm}. 
\end{rmk}

\bigskip

\subsection{Examples for $\ell =9$}

In this section, we shall see some examples in the case $\ell =9$. Here we consider three cases where the coefficients $a_i$ belong to the sets $\{1,2\}$, $\{1,3\}$ and $\{1,2,4\}$. In the first two cases there are 8 forms each corresponding to the vectors ${\bf a} = (1^{i_1}, b^{i_2})$, where $b = 2,3$ and $1\le i_1\le 8$, with $i_1+i_2 =9$. In the third case, there are 36 cases: 
${\bf a} = (1^{i_1}, 2^{i_2},4^{i_3})$, $i_1+i_2+i_3 =9$,  (8 cases if $i_2=0$, 12 cases if $i_2, i_3\not=0$, $i_2$ is even and 
16 cases if $i_2,i_3\not=0$, $i_2$ is odd). For the first two sets, the theta series corresponding to these vectors belong to the space $M_{9/2}(4b,\chi)$, where $\chi$ is the character $\chi_b^{i_2}$, $b=2,3$. For the third set, the theta series belongs to $M_{9/2}(16, \chi_2^{i_2})$.
For the application of \thmref{thm:main}, we can take $t = D =1$ and we see that the Shimura map ${\mathcal S}_1$ maps $\Theta_{\bf a}(z)$ to the space $M_8(2M)$, $M=2,3,4$, respectively corresponding to the three sets.  Therefore, by constructing  bases for the spaces $M_8(N)$, $N=4,6,8$, we derive the formulas for the number of representations of 
$r_9({\bf a};n^2)$ corresponding to these 52 vectors ${\bf a}$. Below we give the necessary details and state the formulas as corollaries. 
  
\bigskip

\centerline{List of basis for $M_8(N)$, $N=4,6,8$}
 
\begin{center}
\begin{tabular}{|c|c|}
\hline 
Space & Basis elements \\
\hline 
$M_8(4)$ & $E_8(z), E_8(2z), E_8(4z),  \Delta_{8,2}(z), \Delta_{8,2}(2z)$ \\
\hline 
$M_8(6)$ & $E_8(z), E_8(2z), E_8(3z),  E_8(6z), \Delta_{8,2}(z), \Delta_{8,2}(3z), \Delta_{8,3}(z), \Delta_{8,3}(2z), 
\Delta_{8,6}(z)$ \\
\hline 
$M_8(8)$ & $E_8(z), E_8(2z), E_8(4z), E_8(8z),  \Delta_{8,2}(z), \Delta_{8,2}(2z), \Delta_{8,2}(4z), \Delta_{8,8;1}(z),   
\Delta_{8,8;2}(z)$\\
\hline 
\end{tabular}
\end{center}

\noindent 
The newforms appearing in the above table are given below.
\begin{align}
\Delta_{8,2}(z) &= \eta^8(z)\eta^8(2z) ~=~ \sum_{n\ge 1} \tau_{8,2}(n) q^n, \\
\Delta_{8,3}(z) & = \eta^{12}(z)\eta^4(3z) + 18 \eta^{9}(z)\eta^4(3z)\eta^3(9z) + 81 \eta^{6}(z)\eta^4(3z)\eta^6(9z) ~
= \sum_{n=1}^{\infty} \tau_{8,3}(n)q^n,\\
\Delta_{8,6}(z) & =\frac{1}{240} \Big( E_4(z)E_4(6z) - E_4(2z)E_4(3z) \Big) = \sum_{n=1}^{\infty} \tau_{8,6}(n)q^n. 
\end{align}
In the case of level $8$, there are two newforms denoted as $\Delta_{8,8;j}(z)$, $j=1,2$, which are given below. 

\begin{align}
\begin{split}
\Delta_{8,8;1}(z) & = \frac{7}{3231360} E_8(z) - \frac{1687}{3231360} E_8(2z) - \frac{1687}{201960}E_8(4z) + \frac{224}{25245} E_8(8z) + \frac{-1385}{2244} \Delta_{8,2}(z) \\
&\quad  +\frac{7450}{561}\Delta_{8,2}(2z) + \frac{15680}{51}\Delta_{8,2}(4z)+\frac{128}{33} \Delta_{8,2}(z) E_4(8z) - \frac{224}{99}G_{4,8}(z),\\
\end{split}
\end{align}
\begin{align}
\begin{split}
\Delta_{8,8;2}(z) & = - \frac{1}{489600} E_8(z) + \frac{241}{489600} E_8(2z) + \frac{241}{30600}E_8(4z) - \frac{32}{3825} E_8(8z) - \frac{77}{68} \Delta_{8,2}(z)\\
& \quad  - \frac{446}{17}\Delta_{8,2}(2z) - \frac{4928}{17}\Delta_{8,2}(4z) + \frac{32}{15}G_{4,8}(z),\\
\end{split}
\end{align}
where $G_{4,8}(z)$ is defined as
\begin{equation}
G_{4,8}(z) = \frac{1}{240}\Big( E_4(z) E_4(8z) - E_4(2z) E_4(4z)\Big).
\end{equation}

\bigskip

For ${\bf a} = (1^{i_1},b^{i_2})$, $i_1+i_2=9$, $b=2,3$, we use the above bases and express the image of $\Theta_{\bf a}(z)$ under the Shimura map ${\mathcal S}_1$. Explicit Shimura images are given in the following table (Table 7 for $b=2$, and Table 8 for $b=3$). When ${\bf a} = (1^{i_1}, 2^{i_2},4^{i_3})$, $i_1+i_2+i_3 =9$, $i_j\ge 1$ we give the Shimura images for these vectors in Table 9 ($i_3\not=0$, $i_2$ is even, including the case $i_2=0$) and in Table 10 ($i_3\not=0$, $i_2$ is odd). In Table 10, we give only the linear combination coefficients $\lambda_{8,8;j}$, $1\le j\le 9$. These are the constants when the Shimura image is written as a linear combination of basis elements listed above for the space $M_8(8)$.\\

\centerline{{\bf Table 7}}

\begin{center}
\begin{tabular}{ |c | c| }
\hline
${\bf a}$ & ${\bf \mathcal{S}}_1({\bf \Theta_a})(z)$  \\ 
\hline
$(1^7,2^2)$ & $ \frac{13}{816} E_8(z)  - \frac{23}{510} E_8(2z) + \frac{108}{17} \Delta_{8,2}(z)$\\  
\hline
$(1^5,2^4)$ & $ \frac{11}{1360} E_8(z)  - \frac{19}{510} E_8(2z) + \frac{104}{17} \Delta_{8,2}(z)$ \\  
\hline
$(1^3,2^6)$ &$ \frac{1}{240} E_8(z)  - \frac{1}{30} E_8(2z) + 4 \Delta_{8,2}(z)$ \\  
\hline
$(1,2^8)$ &  $ \frac{3}{1360} E_8(z)  - \frac{8}{255} E_8(2z) + \frac{16}{17} \Delta_{8,2}(z)$ \\ 
\hline
$(1^2,2^7)$ & $ \frac{11}{4080} E_8(z) - \frac{33}{1360} E_8(2z) +  \frac{1408}{255}  E_8(4z) +  \frac{46}{17} \Delta_{8,2}(z)+ \frac{176}{17} \Delta_{8,2}(2z)$ \\  
\hline
$(1^4,2^5)$ & $ \frac{11}{2040} E_8(z) - \frac{11}{408} E_8(2z) +  \frac{1408}{255}  E_8(4z) +  \frac{92}{17} \Delta_{8,2}(z)+ \frac{448}{17} \Delta_{8,2}(2z)$ \\  
\hline
$(1^6,2^3)$ & $ \frac{11}{1020} E_8(z) - \frac{11}{340} E_8(2z) +  \frac{1408}{255}  E_8(4z) +  \frac{116}{17} \Delta_{8,2}(z)+ \frac{448}{17} \Delta_{8,2}(2z)$ \\  
\hline
$(1^8,2)$ & $ \frac{11}{510} E_8(z) - \frac{11}{255} E_8(2z) +  \frac{1408}{255}  E_8(4z) +  \frac{96}{17} \Delta_{8,2}(z)+ \frac{1536}{17} \Delta_{8,2}(2z)$ \\  
\hline
\end{tabular}
\end{center}

\bigskip

\bigskip

\centerline{{\bf Table 8}}
\smallskip
\begin{center}
\begin{tabular}{ |c | c| }
\hline
${\bf a}$ & ${\bf \mathcal{S}}_1({\bf \Theta_a})(z)$  \\ 
\hline
$(1^7,3^2)$ &$ \frac{91}{9840} E_8(z)  + \frac{2457}{3280} E_8(3z) + \frac{392}{41} \Delta_{8,3}(z)$ \\  
\hline
$(1^5,3^4)$ & $  \frac{137}{33456} E_8(z) -  \frac{16}{2091} E_8(2z) - \frac{2877}{3280} E_8(3z)  + \frac{336}{205} E_8(6z) + \frac{224}{51}\Delta_{8,2}(z)$ \\  
& $ +288 \Delta_{8,2}(3z) + \frac{448}{123} \Delta_{8,3}(z) - \frac{8192}{123} \Delta_{8,3}(2z)$ \\
\hline
$(1^3,3^6)$ &$ \frac{7}{9840} E_8(z)  + \frac{497}{656} E_8(3z) + \frac{232}{41} \Delta_{8,2}(z)$ \\  
\hline
$(1,3^8)$ &  $  \frac{137}{167280} E_8(z) -  \frac{16}{10455} E_8(2z) - \frac{146179}{167280} E_8(3z)  + \frac{17072}{10455} E_8(6z) + \frac{224}{255}\Delta_{8,2}(z)$ \\  
& $ + \frac{1568}{85} \Delta_{8,2}(3z) + \frac{448}{615} \Delta_{8,3}(z) - \frac{8192}{615} \Delta_{8,3}(2z)$ \\
\hline
$(1^2,3^7)$ &  $  \frac{587}{418200} E_8(z) -  \frac{33569}{313650} E_8(2z) - \frac{1106767}{418200} E_8(3z)  + \frac{11106157}{156825} E_8(6z) - \frac{784}{17}\Delta_{8,2}(z)$ \\  
& $ - \frac{51072}{17} \Delta_{8,2}(3z) + \frac{15736}{369} \Delta_{8,3}(z) - \frac{316352}{369} \Delta_{8,3}(2z) +  \frac{4144}{45} \Delta_{8,6}(z)$ \\
\hline
$(1^4,3^5)$ &  $  \frac{16567}{5645700} E_8(z) -  \frac{315263}{1411425} E_8(2z) - \frac{563591}{332100} E_8(3z)  + \frac{9684934}{83025} E_8(6z) - \frac{55424}{459}\Delta_{8,2}(z)$ \\  
& $ - \frac{73952}{9} \Delta_{8,2}(3z) - \frac{130880}{1107} \Delta_{8,3}(z) - \frac{2497024}{1107} \Delta_{8,3}(2z) +  \frac{33152}{135} \Delta_{8,6}(z)$ \\
\hline
$(1^6,3^3)$ & $  \frac{10423}{1254600} E_8(z) -  \frac{147467}{313650} E_8(2z) - \frac{7181323}{1254600} E_8(3z)  + \frac{24943411}{156825} E_8(6z) - \frac{8848}{51}\Delta_{8,2}(z)$ \\  
& $ - \frac{192128}{17} \Delta_{8,2}(3z) - \frac{19528}{123} \Delta_{8,3}(z) - \frac{395840}{123} \Delta_{8,3}(2z) +  \frac{5104}{15} \Delta_{8,6}(z)$ \\
\hline
$(1^8,3)$ & $  \frac{49459}{2822850} E_8(z) -  \frac{89542}{1411425} E_8(2z) - \frac{3071569}{2822850} E_8(3z)  + \frac{146224612}{1411425} E_8(6z) - \frac{218624}{2295}\Delta_{8,2}(z)$ \\  
& $ - \frac{4683328}{765} \Delta_{8,2}(3z) - \frac{511808}{5535} \Delta_{8,3}(z) - \frac{9875968}{5535} \Delta_{8,3}(2z) +  \frac{26368}{135} \Delta_{8,6}(z)$ \\
\hline
\end{tabular}
\end{center}

\bigskip

\bigskip

\smallskip

\bigskip

\begin{center}

\centerline{\bf Table 9}

\bigskip

\begin{tabular}{ |c | c| }
\hline
${\bf a}$ & ${\bf \mathcal{S}}_1({\bf \Theta_a})(z)$  \\ 

\hline

$(1,4^8)$ & $\frac{1}{4080} E_8(z) + \frac{1}{510} E_8(2z) - \frac{8}{255} E_8(4z) + \frac{32}{17} \Delta_{8,2}(z)+ \frac{288}{17} \Delta_{8,2}(2z)$ \\

\hline

$(1^2,4^7)$&$\frac{1}{2040} E_8(z) - \frac{121}{4080}E_8(2z) + \frac{64}{17} \Delta_{8,2}(z)+ 32 \Delta_{8,2}(2z)$ \\

\hline

$(1^3,4^6)$&$\frac{1}{1360} E_8(z) - \frac{25}{408}E_8(2z) + \frac{8}{255}E_8(4z) + \frac{96}{17} \Delta_{8,2}(z)+ \frac{800}{17}\Delta_{8,2}(2z)$ \\

\hline

$(1^4,4^5)$&$\frac{1}{1020} E_8(z) - \frac{251}{4080}E_8(2z) + \frac{8}{255}E_8(4z) + \frac{128}{17} \Delta_{8,2}(z)+ \frac{1072}{17}\Delta_{8,2}(2z)$ \\

\hline

$(1^5,4^4)$&$\frac{1}{680} E_8(z) + \frac{1}{1360} E_8(2z) - \frac{8}{255}E_8(4z) + \frac{158}{17} \Delta_{8,2}(z)+ \frac{1376}{17}\Delta_{8,2}(2z)$ \\

\hline

$(1^6,4^3)$&$\frac{1}{340} E_8(z) + \frac{509}{4080} E_8(2z) - \frac{8}{51}E_8(4z) + \frac{180}{17} \Delta_{8,2}(z)+ \frac{1712}{17}\Delta_{8,2}(2z)$ \\

\hline

$(1^7,4^2)$&$\frac{7}{1020} E_8(z) + \frac{67}{272} E_8(2z) - \frac{24}{85}E_8(4z) + \frac{182}{17} \Delta_{8,2}(z)+ \frac{2048}{17}\Delta_{8,2}(2z)$ \\

\hline

$(1^8,4)$&$\frac{4}{255} E_8(z) + \frac{19}{80} E_8(2z) - \frac{24}{85}E_8(4z) + \frac{144}{17} \Delta_{8,2}(z)+ \frac{2048}{17}\Delta_{8,2}(2z)$ \\

\hline

$(1,2^2,4^6)$&$\frac{1}{4080} E_8(z) + \frac{1}{510} E_8(2z) - \frac{8}{255}E_8(4z) + \frac{32}{17} \Delta_{8,2}(z)+ \frac{288}{17}\Delta_{8,2}(2z)$ \\

\hline

$(1,2^4,4^4)$&$\frac{1}{2040} E_8(z) + \frac{7}{4080} E_8(2z) - \frac{8}{255}E_8(4z) + \frac{30}{17} \Delta_{8,2}(z)+ \frac{288}{17}\Delta_{8,2}(2z)$ \\

\hline

$(1,2^6,4^2)$&$\frac{1}{1020} E_8(z) + \frac{1}{816} E_8(2z) - \frac{8}{255}E_8(4z) + \frac{26}{17} \Delta_{8,2}(z)+ \frac{288}{17}\Delta_{8,2}(2z)$ \\

\hline

$(1^2,2^2,4^5)$&$\frac{1}{2040} E_8(z) + \frac{7}{4080} E_8(2z) - \frac{8}{255}E_8(4z) + \frac{64}{17} \Delta_{8,2}(z)+ \frac{560}{17}\Delta_{8,2}(2z)$ \\

\hline

$(1^2,2^4,4^3)$&$\frac{1}{1020} E_8(z) + \frac{1}{816} E_8(2z) - \frac{8}{255}E_8(4z) + \frac{60}{17} \Delta_{8,2}(z)+ \frac{560}{17}\Delta_{8,2}(2z)$ \\

\hline

$(1^2,2^6,4^1)$&$\frac{1}{510} E_8(z) + \frac{1}{4080} E_8(2z) - \frac{8}{255}E_8(4z) + \frac{52}{17} \Delta_{8,2}(z)+ \frac{288}{17}\Delta_{8,2}(2z)$ \\

\hline

$(1^3,2^2,4^4)$&$\frac{1}{1020} E_8(z) + \frac{1}{816} E_8(2z) - \frac{8}{255}E_8(4z) + \frac{94}{17} \Delta_{8,2}(z)+ \frac{832}{17}\Delta_{8,2}(2z)$ \\

\hline

$(1^3,2^4,4^2)$&$\frac{1}{510} E_8(z) + \frac{1}{4080} E_8(2z) - \frac{8}{255}E_8(4z) + \frac{86}{17} \Delta_{8,2}(z)+ \frac{832}{17}\Delta_{8,2}(2z)$ \\

\hline

$(1^4,2^2,4^3)$&$\frac{1}{510} E_8(z) + \frac{1}{4080} E_8(2z) - \frac{8}{255}E_8(4z) + \frac{120}{17} \Delta_{8,2}(z)+ \frac{1104}{17}\Delta_{8,2}(2z)$ \\

\hline

$(1^4,2^4,4)$&$\frac{1}{255} E_8(z) - \frac{7}{4080} E_8(2z) - \frac{8}{255}E_8(4z) + \frac{104}{17} \Delta_{8,2}(z)+ \frac{832}{17}\Delta_{8,2}(2z)$ \\

\hline

$(1^5,2^2,4^2)$&$\frac{1}{255} E_8(z) - \frac{7}{4080} E_8(2z) - \frac{8}{255}E_8(4z) + \frac{138}{17} \Delta_{8,2}(z)+ \frac{1376}{17}\Delta_{8,2}(2z)$ \\

\hline

$(1^6,2^2,4)$&$\frac{2}{255} E_8(z) - \frac{23}{4080} E_8(2z) - \frac{8}{255}E_8(4z) + \frac{140}{17} \Delta_{8,2}(z)+ \frac{1376}{17}\Delta_{8,2}(2z)$ \\

\hline

\end{tabular}

\end{center}

\smallskip
\bigskip
\bigskip

\centerline{\bf Table 10}
\smallskip

\begin{center}

\begin{tabular}{|c|c|c|c|c|c|c|c|c|c|}

\hline
$ {\bf a} $&$\lambda_{8,8;1}$&$\lambda_{8,8;2}$&$\lambda_{8,8;3}$&$\lambda_{8,8;4}$&$\lambda_{8,8;5}$&$\lambda_{8,8;6}$&$\lambda_{8,8;7}$&$\lambda_{8,8;8}$&$\lambda_{8,8;9}$\\
\hline
&&&&&&&&&\\
$(1,2,4^7)$&$\frac{11}{65280}$&$\frac{-11}{65280}$&$\frac{-11}{510}$&$\frac{1408}{255}$&$\frac{125}{136}$&13&$\frac{1536}{17}$&$\frac{7}{8}$&$\frac{1}{8}$\\
\hline
&&&&&&&&&\\
$(1,2^3,4^5)$&$\frac{11}{32640}$&$\frac{-11}{32640}$&$\frac{-11}{510}$&$\frac{1408}{255}$&$\frac{57}{68}$&10&$\frac{3712}{17}$&1& 0 \\
\hline
&&&&&&&&&\\
$(1,2^5,4^3)$&$\frac{11}{16320}$&$\frac{-11}{16320}$&$\frac{-11}{510}$&$\frac{1408}{255}$&$\frac{23}{34}$&12&$\frac{5888}{17}$&1& 0 \\
\hline
&&&&&&&&&\\
$(1,2^7,4)$&$\frac{11}{8160}$&$\frac{-11}{8160}$&$\frac{-11}{510}$&$\frac{1408}{255}$&$\frac{23}{17}$&16&$\frac{5888}{17}$&0& 0 \\
\hline
&&&&&&&&&\\
$(1^2,2,4^6)$&$\frac{11}{32640}$&$\frac{-11}{32640}$&$\frac{-11}{510}$&$\frac{1408}{255}$&$\frac{125}{68}$&18&$\frac{3712}{17}$&$\frac{7}{4}$&$\frac{1}{4}$\\
\hline
&&&&&&&&&\\
$(1^2,2^3,4^4)$&$\frac{11}{16320}$&$\frac{-11}{16320}$&$\frac{-11}{510}$&$\frac{1408}{255}$&$\frac{57}{34}$&20&$\frac{5888}{17}$&2& 0 \\
\hline
&&&&&&&&&\\
$(1^2,2^5,4^2)$&$\frac{11}{8160}$&$\frac{-11}{8160}$&$\frac{-11}{510}$&$\frac{1408}{255}$&$\frac{23}{17}$&16&$\frac{5888}{17}$&2& 0 \\
\hline
\end{tabular}
\end{center}

\smallskip

\bigskip

\begin{center}

Table 10 contd. \\
\smallskip

\begin{tabular}{|c|c|c|c|c|c|c|c|c|c|}
\hline
&&&&&&&&&\\
$(1^3,2,4^5)$&$\frac{11}{16320}$&$\frac{-11}{16320}$&$\frac{-11}{510}$&$\frac{1408}{255}$&$\frac{91}{34}$&20&$\frac{8064}{17}$&$\frac{11}{4}$&$\frac{1}{4}$\\
\hline
&&&&&&&&&\\
$(1^3,2^3,4^3)$&$\frac{11}{8160}$&$\frac{-11}{8160}$&$\frac{-11}{510}$&$\frac{1408}{255}$&$\frac{40}{17}$&32&$\frac{10240}{17}$&3& 0 \\
\hline
&&&&&&&&&\\
$(1^3,2^5,4)$&$\frac{11}{4080}$&$\frac{-11}{4080}$&$\frac{-11}{510}$&$\frac{1408}{255}$&$\frac{46}{17}$&24&$\frac{5888}{17}$&2& 0 \\
\hline
&&&&&&&&&\\
$(1^4,2,4^4)$&$\frac{11}{8160}$&$\frac{-11}{8160}$&$\frac{-11}{510}$&$\frac{1408}{255}$&$\frac{57}{17}$&24&$\frac{14592}{17}$&4& 0 \\
\hline
&&&&&&&&&\\
$(1^4,2^3,4^2)$&$\frac{11}{4080}$&$\frac{-11}{4080}$&$\frac{-11}{510}$&$\frac{1408}{255}$&$\frac{46}{17}$&40&$\frac{14592}{17}$&4& 0 \\
\hline
&&&&&&&&&\\
$(1^5,2,4^3)$&$\frac{11}{4080}$&$\frac{-11}{4080}$&$\frac{-11}{510}$&$\frac{1408}{255}$&$\frac{63}{17}$&40&$\frac{23296}{17}$&$\frac{11}{2}$&$\frac{-1}{2}$\\
\hline
&&&&&&&&&\\
$(1^5,2^3,4)$&$\frac{11}{2040}$&$\frac{-11}{2040}$&$\frac{-11}{510}$&$\frac{1408}{255}$&$\frac{58}{17}$&40&$\frac{14592}{17}$&4& 0 \\
\hline
&&&&&&&&&\\
$(1^6,2,4^2)$&$\frac{11}{2040}$&$\frac{-11}{2040}$&$\frac{-11}{510}$&$\frac{1408}{255}$&$\frac{58}{17}$&72&$\frac{32000}{17}$&7& -1 \\
\hline
&&&&&&&&&\\
$(1^7,2,4)$&$\frac{11}{1020}$&$\frac{-11}{1020}$&$\frac{-11}{510}$&$\frac{1408}{255}$&$\frac{48}{17}$&104&$\frac{32000}{17}$&7& -1 \\
\hline 
\end{tabular}
\end{center}
\smallskip

\bigskip

\smallskip

\noindent {\bf Simplifications and some congruence relations}: 
We now present simplified versions of some of the formulas (4 from Table 7, 2 each from Tables 8 and 9). These are corresponding to the vectors $(1^7,2^2), (1^5,2^4), (1^3, 2^6), (1, 2^8)$ from Table 7, $(1^7,3), (1^3, 3^6)$ from Table 8 and $(1,4^8), (1^2,4^7)$ from Table 9. \\

\noindent  
For an odd positive integer $m$, let 
\begin{align}
C_5(m) &= \tau_{8,2}(m) \prod_{p\ge 3}\left(1 - p^3 \frac{\tau_{8,2}(m/p)}{\tau_{8,2}(m)}\right), \label{c3}\\
C_6(m) & = \tau_{8,3}(m) \prod_{p\ge 5}\left(1 - p^3 \frac{\tau_{8,3}(m/p)}{\tau_{8,3}(m)}\right). \label{c4}
\end{align}

\noindent We use the same notation for the prime factorisation of a natural number $n$, with $\lambda_p$ denoting the highest power of $p$ dividing $n$ and $m$ is the odd part of $n$. The 8 formulas are listed in the following corollary.

\begin{cor}\label{9var}
For $n = 2^{\lambda_2}m$, $m\ge 1$, we have 
\begin{align}
r_9(1^7,2^2; n^2) & = \frac{2}{17}\left(\frac{1017\times 2^{7\lambda_2+3} +119}{2^7-1}\right) ~s_7(m) + 
\frac{108}{17} \tau_{8,2}(2^{\lambda_2}) ~ C_5(m),\label{64}\\
r_9(1^5,2^4; n^2) & = \frac{2}{17}\left(\frac{509\times 2^{7\lambda_2+3} +119}{2^7-1}\right) ~s_7(m) + 
\frac{104}{17} \tau_{8,2}(2^{\lambda_2})  C_5(m), \label{65}\\
r_9(1^3,2^6; n^2) & = 2 \left(\frac{15\times 2^{7\lambda_2+3} +7}{2^7-1}\right) ~s_7(m) + 2 \tau_{8,2}(2^{\lambda_2}) 
C_5(m),\label{66}
\end{align}
\smallskip
\begin{align}
r_9(1,2^8; n^2) & = \frac{2}{17}\left(\frac{8\times 2^{7\lambda_2+7} +119}{2^7-1}\right) ~s_7(m)  + \frac{8}{17} 
\tau_{8,2}(2^{\lambda_2})  C_5(m),\label{67}\\
r_9(1,4^8;n^2) & = \frac{2}{17} \left(\frac{67\times 2^{7\lambda_2+1}+247}{2^7-1}\right) s_7(m) + \frac{32}{17} \big(\tau_{8,2}(2^{\lambda_2}) + 9 \tau_{8,2}(2^{\lambda_2-1})\big) C_5(m),\label{70}\\
r_9(1^2,4^7;n^2) & =  \frac{2}{17} \left(\frac{135\times 2^{7\lambda_2}+119}{2^7-1}\right) s_7(m) + \big(\frac{64}{17} 
\tau_{8,2}(2^{\lambda_2}) +  32 \tau_{8,2}(2^{\lambda_2-1})\big) C_5(m),\label{71}
\end{align}
where $s_7(m)$ is defined by \eqref{s7} and $C_5(m)$ is defined by \eqref{c3}. 
When $n= 2^{\lambda_2} 3^{\lambda_3} m$, with $\gcd(m,6)=1$, we have 
\begin{align}
r_9(1^7,3^2;n^2) & = \frac{364}{41} \sigma_7(2^{\lambda_2}) \left(\frac{14\times 3^{7\lambda_3+4} -41}{3^7-1}\right) 
s_7(m) + \frac{392}{41} \tau_{8,3}(2^{\lambda_2}3^{\lambda_3}) C_6(m),\label{68}\\
r_9(1^3,3^6;n^2) & = \frac{14}{41} \sigma_7(2^{\lambda_2}) \left(\frac{1084\times 3^{7\lambda_3+1} -1066}{3^7-1}\right) s_7(m)  + \frac{232}{41}\tau_{8,3}(2^{\lambda_2}3^{\lambda_3}) C_6(m),\label{69}
\end{align}
where $C_6(m)$ is defined by \eqref{c4}.
\end{cor}

\bigskip

From the above formulas, one can derive certain congruences for the Fourier coefficients of the two newforms of weight $8$ 
and levels $2$ and $3$. \\

\begin{cor}\label{cong}
Let $\lambda \ge 1$ be an integer. Then for an odd prime $p$, we have 
\begin{align}
\tau_{8,2}(p^\lambda) &\equiv \sigma_7(p^\lambda)\pmod{17}. \label{cong1}
\end{align}
For odd primes $p\ge 5$, we have 
\begin{align}
\tau_{8,3}(p^\lambda) &\equiv \sigma_7(p^\lambda)\pmod{41}. \label{cong2}
\end{align}
In particular, for any integer $\lambda \ge 1$, we have 
\begin{align}
\tau_{8,2}(17^\lambda) & \equiv 1\pmod{17}, \label{cong3}\\
\tau_{8,3}(41^\lambda) & \equiv 1\pmod{41}. \label{cong4}
\end{align}
\end{cor}

\smallskip

\begin{rmk}
The cases in Table 9 correspond to forms of weight $9/2$ on $\Gamma_0(16)$ with trivial character and so by the Shimura map these forms are mapped to the space $M_8(8)$ and we used a basis of this space to express the images. However, by looking at the image functions, all of them are in the subspace $M_8(4)$. This may be due to some properties, which can be explored. 
\end{rmk}

\bigskip

\section{Concluding remarks}
 
\smallskip
We have mentioned that in some cases we could not get the required formulas for the number of representations of $n^2$ and instead we could get only the formulas for $4n^2$. For example the case when $\ell \equiv3 \pmod{4}$ and $N_{\bf a}$ is even, we only have the mapping property of the Shimura-Kohnen map.  However, when we consider a square-free integer $t\equiv 3\pmod{4}$, then by taking $D =-t$, which is a negative fundamental discriminant, the application of the Shimura-Kohnen map gives rise to formulas for $r_\ell({\bf a}; |t|n^2)$. So, instead of formulas for $n^2$ we get formulas for $|t|n^2$, where $t$ is a square-free integer and $-t$ is a fundamental discriminant. In the case of quadratic forms of odd number of squares, only representation numbers of squares are known so far and formulas for non-square integers are not yet proved. Though there are some restrictions, the method of Shimura correspondence gives formulas for a wider class. One also gets formulas in the case when $\ell \equiv 1\pmod{4}$ (in this case we can take $t$ to be positive square-free which is congruent to 1 modulo 4). 
Towards this direction, for the case $\ell=7$, we present some formulas for the number of representations of $3n^2$ for certain coefficient vectors ${\bf a}$. \\

We consider the coefficients in the sets $\{1,2\}$, $\{1,3\}$ and $\{1,2,4\}$. By using the basis respectively for the spaces $M_6(N)$, $N=4,6,8$, we obtain the following formulas. The basis elements of $M_6(4)$ and $M_6(8)$ were given in \S 4.2. Here we give a basis for $M_6(6)$, consisting the following $7$ modular forms: $E_6(z)$, $E_6(2z)$, $E_6(3z)$, $E_6(6z)$, $\Delta_{6,3}(z)$, $\Delta_{6,3}(2z)$, $\Delta_{6,6}(z)$. To derive our formulas, we denote these bases by $f_{N,j}(z)$, $1\le j\le \dim_{\mathbb C}(M_6(N))$, where $N = 4,6,8$ and the dimensions of the spaces are $4,7,7$ respectively. So, by the application of \thmref{thm:main}, we get the following expression (taking $D=-3$, $\ell =7$).\\

\begin{equation}
{\mathcal S}_{-3}(\Theta_{\bf a}(z)) = \sum_j \lambda_{{\bf a}, N, j} f_{N,j}(z),
\end{equation}
where the constants $\lambda_{{\bf a}, N, j}$ depend on the vector ${\bf a}$. By comparing the $n$-th Fourier coefficients and taking the M\"obius inversion, we get the following formula: 
\begin{equation}
r_7({\bf a}; 3n^2) = \sum_{d\vert n \atop{\gcd(d,2N_{\bf a})=1}} \mu(d) \psi_{\bf a}(d) d^2 \sum_{j} \lambda_{{\bf a},N,j} a_{f_{N,j}}(n/d),
\end{equation}
where $j$ runs from $1$ to the dimension of the respective space, $a_{f_{N,j}}(n)$ denotes the $n$-th Fourier coefficients of the basis elements. 

\bigskip

In the following tables, we give the explicit values of the constants $\lambda_{{\bf a},N,j}$, where ${\bf a}$ belongs to the three sets $\{1,2\}$, $\{1,3\}$, $\{1,2,4\}$. \\

\vfill

\newpage 

\begin{center}

{\tiny 
\textbf{Table 11.1 for the coefficients $ \lambda_{{\bf a}, 4, j}$, $1\le j\le 4$, ${\bf a} \in \{1,2\}$}
\vspace{0.2cm}

\begin{tabular}{|c|c|c|c|c|c|c|}
\hline
$\rightarrow j $&$1$&$2$&$3$&$4$\\
$\downarrow {\bf a}$&&&&\\
\hline
&&&&\\
$1^5,2^2$&$\frac{-5}{21}$&$\frac{-20}{63}$&0&0\\
&&&&\\
$1^3,2^4$&$\frac{-1}{9}$&$\frac{-4}{9}$&0&0\\
&&&&\\
$1,2^6$&$\frac{-1}{21}$&$\frac{-32}{63}$&0&0\\
&&&&\\
\hline
\end{tabular}
\qquad 
\begin{tabular}{|c|c|c|c|c|c|c|}
\hline
$\rightarrow j $&$1$&$2$&$3$&$4$\\
$\downarrow {\bf a}$&&&&\\
\hline
&&&&\\
$1^6,2$&$\frac{-23}{63}$&$\frac{46}{63}$&$\frac{-1472}{63}$& 0\\
&&&&\\
$1^4,2^3$&$\frac{-23}{126}$&$\frac{23}{42}$&$\frac{-1472}{63}$& -12\\
&&&&\\
$1^2,2^5$&$\frac{-23}{252}$&$\frac{115}{252}$&$\frac{-1472}{63}$& -6\\
&&&&\\
\hline
\end{tabular}
}
\end{center}

\begin{center}
{\tiny 
\textbf{Table 11.2 for the coefficients $ \lambda_{{\bf a}, 6, j}$, $1\le j\le 7$, ${\bf a} \in \{1,3\}$}
\vspace{0.2cm}

\begin{tabular}{|c|c|c|c|c|c|c|c|}
\hline
$\rightarrow j $&$1$&$2$&$3$&$4$&$5$&$6$&$7$\\
$\downarrow {\bf a}$&&&&&&&\\
\hline
&&&&&&&\\
$1^5,3^2 $&$\frac{-1769}{11466}$&$\frac{1952}{5733}$&$\frac{783}{1274}$&$\frac{-864}{637}$&$\frac{200}{13}$&$\frac{1280}{13}$&$\frac{-64}{7}$ \\
&&&&&&&\\
$1^3,3^4 $&$\frac{-50}{819}$&0&$\frac{-45}{91}$&0&$\frac{-192}{13}$&0&0 \\
&&&&&&&\\
$1,3^6 $&$\frac{-29}{1638}$&$\frac{32}{819}$&$\frac{87}{182}$&$\frac{-96}{91}$&$\frac{40}{13}$&$\frac{256}{13}$&0 \\
&&&&&&&\\
\hline
\end{tabular}
\quad 
\begin{tabular}{|c|c|c|c|c|c|c|c|}
\hline
$\rightarrow j $&$1$&$2$&$3$&$4$&$5$&$6$&$7$\\
$\downarrow {\bf a}$&&&&&&&\\
\hline
&&&&&&&\\
$1^6,3 $&$\frac{-113}{364}$&$\frac{226}{91}$&$\frac{243}{364}$&$\frac{-486}{91}$&$\frac{72}{13}$&$\frac{-576}{13}$&0 \\
&&&&&&&\\
$1^4,3^3 $&$\frac{-49}{468}$&$\frac{-98}{117}$&$\frac{-9}{52}$&$\frac{-18}{13}$&$\frac{-192}{13}$&$\frac{-1536}{13}$&0 \\
&&&&&&&\\
$1^2,3^5 $&$\frac{-5}{364}$&$\frac{10}{91}$&$\frac{135}{364}$&$\frac{-270}{91}$&$\frac{40}{13}$&$\frac{-320}{13}$&0 \\ 
&&&&&&&\\
\hline
\end{tabular}
}
\end{center}

\begin{center}
{\tiny 
\textbf{Table 11.3 for the coefficients $ \lambda_{{\bf a}, 8, j}$, $1\le j\le 7$, ${\bf a} \in \{1,2,4\}$}
\vspace{0.2cm}

\begin{tabular}{|c|c|c|c|c|c|c|c|c|c|}
\hline
$\rightarrow j $&$1$&$2$&$3$&$4$&$5$&$6$&$7$\\
$\downarrow {\bf a}$&&&&&&&\\
\hline
&&&&&&&\\
$1^6,4 $&$\frac{-20}{63}$&$\frac{145}{63}$&$\frac{-160}{63}$&0& 0 & 0 & 0 \\
&&&&&&&\\
$1^5,4^2 $&$\frac{-10}{63}$&$\frac{15}{7}$&$\frac{-160}{63}$&0& 0 & 0 & 0 \\
&&&&&&&\\
$1^4,4^3 $&$\frac{-4}{63}$&$\frac{65}{63}$&$\frac{-32}{21}$&0& 0 & 0 & 0 \\
&&&&&&&\\
$1^3,4^4 $&$\frac{-1}{63}$&$\frac{-2}{63}$&$\frac{-32}{63}$&0& 0 & 0 & 0 \\
&&&&&&&\\
$1^2,4^5 $& 0&$\frac{-5}{9}$& 0& 0 & 0 & 0 & 0 \\
&&&&&&&\\
$1, 4^6 $& 0&$\frac{-5}{9}$& 0& 0 & 0 & 0 & 0 \\
&&&&&&&\\
\hline
\end{tabular}
\quad
\begin{tabular}{|c|c|c|c|c|c|c|c|}
\hline
$\rightarrow j $&$1$&$2$&$3$&$4$&$5$&$6$&$7$\\
$\downarrow {\bf a}$&&&&&&&\\
\hline
&&&&&&&\\
$1^4,2^2,4 $&$\frac{-8}{63}$&$\frac{5}{63}$&$\frac{-32}{63}$&0& 0 & 0 & 0 \\
&&&&&&&\\
$1^3,2^2,4^2 $&$\frac{-4}{63}$&$\frac{1}{63}$&$\frac{-32}{63}$&0& 0 & 0 & 0 \\
&&&&&&&\\
$1^2,2^4,4 $&$\frac{-4}{63}$&$\frac{1}{63}$&$\frac{-32}{63}$&0& 0 & 0 & 0 \\
&&&&&&&\\
$1^2,2^2,4^3 $&$\frac{-2}{63}$&$\frac{-1}{63}$&$\frac{-32}{63}$&0& 0 & 0 & 0 \\
&&&&&&&\\
$1,2^4,4^2 $&$\frac{-2}{63}$&$\frac{-1}{63}$&$\frac{-32}{63}$&0& 0 & 0 & 0 \\
&&&&&&&\\
$1,2^2,4^4 $&$\frac{-1}{63}$&$\frac{-2}{63}$&$\frac{-32}{63}$&0& 0 & 0 & 0 \\
&&&&&&&\\
\hline
\end{tabular}
}
\end{center}


\begin{center}
{\tiny 
\begin{tabular}{|c|c|c|c|c|c|c|c|}
\hline
$\rightarrow j $&$1$&$2$&$3$&$4$&$5$&$6$&$7$\\
$\downarrow {\bf a}$&&&&&&&\\
\hline
&&&&&&&\\
$1^5,2,4 $&$\frac{-23}{126}$&$\frac{23}{126}$&$\frac{23}{63}$&$\frac{-1472}{63}$&0&120&8 \\
&&&&&&&\\
$1^4,2,4^2 $&$\frac{-23}{252}$&$\frac{23}{252}$&$\frac{23}{63}$&$\frac{-1472}{63}$&-6&120&8 \\
&&&&&&&\\
$1^3,2^3,4 $&$\frac{-23}{252}$&$\frac{23}{252}$&$\frac{23}{63}$&$\frac{-1472}{63}$&-6&24&4 \\
&&&&&&&\\
$1^3,2,4^3 $&$\frac{-23}{504}$&$\frac{23}{504}$&$\frac{23}{63}$&$\frac{-1472}{63}$&-9&72&6 \\
&&&&&&&\\
$1^2,2^3,4^2 $&$\frac{-23}{504}$&$\frac{23}{504}$&$\frac{23}{63}$&$\frac{-1472}{63}$&-3&24&4 \\
&&&&&&&\\
$1^2,2,4^4 $&$\frac{-23}{1008}$&$\frac{23}{1008}$&$\frac{23}{63}$&$\frac{-1472}{63}$&$\frac{-15}{2}$&24&4 \\
&&&&&&&\\
$1,2^5,4 $&$\frac{-23}{504}$&$\frac{23}{504}$&$\frac{23}{63}$&$\frac{-1472}{63}$&-3&24&0 \\
&&&&&&&\\
$1,2^3,4^3 $&$\frac{-23}{1008}$&$\frac{23}{1008}$&$\frac{23}{63}$&$\frac{-1472}{63}$&$\frac{-3}{2}$&24&2 \\
&&&&&&&\\
$1,2,4^5 $&$\frac{-23}{2016}$&$\frac{23}{2016}$&$\frac{23}{63}$&$\frac{-1472}{63}$&$\frac{-15}{4}$&0&2 \\
&&&&&&&\\     	
\hline
\end{tabular}
}
\end{center}

\bigskip

We now give some of the formulas from the above data in explicit form. 

\begin{cor}
For a natural number $n = 2^{\lambda_2}m$, with $m$ odd, we have the following formulas for the number of representations. 
\begin{align}
r_7(1^6,2;3n^2) & = 184 \left|\frac{2^{5\lambda_2+5} -63}{2^5-1}\right| s_5(m),\\
r_7(1^5,2^2; 3n^2) &=  40 \left(\frac{25\times 2^{5\lambda_2+2} - 7}{2^5-1}\right) s_6(m),\\
r_7(1^4,2^3;3n^2) & = \begin{cases} 92 ~s_5(m) - 12 C_3(m) &  {\rm ~if~} \lambda_2 =0, \\
92\left(\frac{33 \times 2^{5\lambda_2} -126}{2^5-1}\right) s_5(m) &  {\rm ~if~} \lambda_2 \ge 1, \\
\end{cases}\\
r_7(1^3,2^4;3n^2) & = 56 \left(\frac{9 \times 2^{5\lambda_2+2} -5}{2^5-1}\right) s_6(m),\\
r_7(1^2,2^5,3n^2) & = \begin{cases} 46 ~s_5(m) - 6 C_3(m) &  {\rm ~if~} \lambda_2 =0, \\
46\left(\frac{35 \times 2^{5\lambda_2} -252}{2^5-1}\right) s_5(m) &  {\rm ~if~} \lambda_2 \ge 1, \\
\end{cases}\\
r_7(1,2^6;3n^2) & = 8 \left(\frac{2^{5\lambda_2 +7} -35}{2^5-1}\right) s_6(m),\\
r_7(1^2,2^2,4^3; 3n^2) & = r_7(1,2^4,4^2;3n^2) ~=~ \begin{cases} 16~ s_6(m)  &  {\rm ~if~} \lambda_2 =0, \\
8\left(\frac{66 \times 2^{5\lambda_2} -35}{2^5-1}\right) s_6(m) &  {\rm ~if~} \lambda_2 \ge 1, \\
\end{cases}\\
r_7(1^3,2^2,4^2; 3n^2) & = r_7(1^2,2^4,4;3n^2) ~=~ \begin{cases} 32 ~s_6(m)  &  {\rm ~if~} \lambda_2 =0, \\
8\left(\frac{2^{5\lambda_2+7} -35}{2^5-1}\right) s_6(m) &  {\rm ~if~} \lambda_2 \ge 1, \\
\end{cases}\\
r_7(1^4,2^2,4; 3n^2) & =  \begin{cases} 64 ~s_6(m)  &  {\rm ~if~} \lambda_2 =0, \\
2072 ~\sigma_5(2^{5\lambda_2-1}) s_6(m)  &{\rm ~if~} \lambda_2 \ge 1, \\
\end{cases}\\
r_7(1^6,4 ; 3n^2) & = \begin{cases} 160~ s_6(m)  &  {\rm ~if~} \lambda_2 =0, \\
40 \left( \frac{25 \times 2^{5\lambda_2+2} -7}{2^5-1} \right) s_6(m) & {\rm ~if~} \lambda_2\ge 1,\\
\end{cases}\\
r_7(1^5,4^2 ; 3n^2)  & = \begin{cases} 80~ s_6(m)  &  {\rm ~if~} \lambda_2 =0, \\
40 \left( \frac{19 \times 2^{5\lambda_2+1} -7}{2^5-1} \right) s_6(m) & {\rm ~if~} \lambda_2\ge 1,\\
\end{cases}\\
r_7(1^4,4^3 ; 3n^2)  & = \begin{cases} 32~ s_6(m)  &  {\rm ~if~} \lambda_2 =0, \\
8 \left( \frac{33 \times 2^{5\lambda_2+1} -35}{2^5-1} \right) s_6(m) & {\rm ~if~} \lambda_2\ge 1,\\
\end{cases}
\end{align}
\begin{align}
r_7(1^3,4^4 ; 3n^2)  & =   r_7(1^2,4^5 ; 3n^2) =  r_7(1,4^6 ; 3n^2) = r_7(1,2^2,4^4; 3n^2) \\
& = 280 ~\sigma_5(2^{\lambda_2-1})  s_6(m), {\rm ~if~} \lambda_2\ge 1,\label{93}\\
r_7(1^3, 4^4; 3m^2) & = 35 ~r_7(1,2^2,4^4; 3m^2) ~=~ 280 ~s_6(m),\\
r_7(1^3,3^4;3n^2) & = \frac{40}{13} \sigma_5(2^{\lambda_2}) \left(\frac{31\times 3^{5\lambda_3+4} - 91}{3^5-1}\right) s_6'(m) - \frac{192}{13} \tau_{6,3}(2^{\lambda_2}3^{\lambda_3}) C_7(m). \label{cong5}
\end{align}
In the last formula, we take $n= 2^{\lambda_2} 3^{\lambda_3}m$, $\gcd(m,6)=1$. The functions $s_5(m)$ and $s_6(m)$ are defined as in \eqref{sm}, \eqref{tm} respectively and the functions $s_6'(m)$ and $C_7(m)$ are defined below.  
\begin{align}
s_6'(m) &= \prod_{p\ge 5} \left(\frac{p^{5\lambda_p+5} -1}{p^5-1} - p^2 \left(\frac{-4}{p}\right) 
\frac{p^{5\lambda_p}-1}{p^5-1}\right), \\
C_7(m) & = \tau_{6,3}(m) \prod_{p\ge 5}\left(1 - p^2 \left(\frac{-4}{p}\right)\frac{\tau_{6,3}(m/p)}{\tau_{6,3}(m)}\right). 
\end{align}
\end{cor}

\bigskip

\begin{rmk}
We remark that numbers of the form $3m^2$, where $m$ is odd are not represented by the forms corresponding to $(1^2,4^5)$ and $(1,4^6)$. That is $r_7(1^2, 4^5;3m^2) = r_7(1, 4^6;3m^2) = 0$, if $m$ is odd. This fact can also be checked by elementary arguments. When $n$ is even these numbers of representations (of $3n^2$) are equal to the number of representations of $3n^2$ by the form corresponding to $(1^3,4^4)$ (see Eq.\eqref{93}). There are many relations among these representation numbers. We will mention some of them here. 
\begin{equation*}
\begin{split}
r_7(1^6,4;3n^2) &= r_7(1^5,2^2; 3n^2) ~({\rm if~} 2\vert n),\\
r_7(1^4;2^3; 3m^2) & =2 ~ r_7(1^2,2^5; 3m^2), \\
r_7(1^4,4^3; 3m^2) & = r_7(1^3,2^2,4^2; 3m^2) ~= r_7(1^2,2^4,4; 3m^2)\\
 &= 2 ~r_7(1^2,2^2,4^3; 3m^2) ~= ~2 r_7(1^2,2^4,4; 3m^2),\\
 r_7(1^6,4;3m^2) & = r_7(1^4,2^2,4;3m^2) \\
 & = 2~r_7(1^4,4^3; 3m^2) ~=~ 2~r_7(1^5,4^2; 3m^2).\\
\end{split}
\end{equation*}

\end{rmk}

\begin{rmk}
There is one more congruence from the last formula in the above corollary. 
It is easy to derive the following congruence using \eqref{cong5}: 
\begin{equation}
\tau_{6,3}(p^\lambda) \equiv \sigma_5(p^\lambda) \pmod{13},
\end{equation}
where $p\not =3$ is a prime and $\lambda \ge 1$ is an integer. 
\end{rmk}

\noindent {\bf Further results}. In our earlier works we determined explicit bases for the spaces of modular forms $M_4(N)$, $N= 4,6,8,12,14,24,28,32,48$ and $M_6(N')$, $N'=4,6,8,12$. Therefore, it is possible to use the method described in our work to include coefficients $a_i$ in larger sets, which we have done and computed explicit description of the images under the Shimura maps in the case when $\ell =5$ (quinary) and $\ell =7$ (septenary). Specifically, in the quinary case we can take $a_i$ in the following sets: $\{1,2,3,4,6,12,24\}$, $\{1,2,4,8,16\}$, $\{1,2,7,14\}$, $\{1,11\}$. In the septenary case  we can take 
$a_i$ in the set: $\{1,2,3,4,6\}$. Detailed list of tables for the linear combination coefficients are available. Since they occupy large space (more than 10 pages), we have not presented in this paper. From these data, one can write down more than 500 explicit formulas for the number of representations of quadratic forms corresponding to these $a_i$'s. 
We will be uploading these details in preprint form for reference. 

\smallskip

As a final remark, we mention that it is difficult to make use of Hurwitz method in the general set up. Also by using this method one gets formulas for only squares. So, at the moment using Shimura map seems to give new and general results. However, our approach is also not giving a complete answer to the original question of obtaining explicit formulas for the number of representations of any natural number.

\bigskip

{\tiny
\noindent {\bf Acknowledgements}:  We have used the open-source mathematics software SAGE (www.sagemath.org) \cite{sage} to perform our calculations. We also used the L-functions and Modular forms database (LMFDB) \cite{lmfdb}. 
Both are acknowledged in the References section. Parts of the work were done when the first and the last authors were at the Harish-Chandra Research Institute, Prayagraj and some parts were carried out during their visits to NISER, Bhubaneswar. They would like to thank both HRI and NISER for their support.}

\end{document}